\documentclass[11pt,review,letterpaper,authoryear]{elsarticle}

\usepackage[english]{babel}
\usepackage[utf8]{inputenc}

\usepackage{geometry}
\usepackage[onehalfspacing]{setspace} 

\usepackage{multicol,float}
\usepackage{color}
\usepackage{multirow}
\usepackage{array}
\usepackage{amsmath}

\usepackage{amsfonts}
\usepackage{amssymb}
\usepackage{graphicx}
\usepackage{lineno}
\usepackage{colortbl}
\usepackage{rotating}
\usepackage{cancel}
\usepackage{efbox}
\usepackage{subcaption}
\usepackage{ulem}
\usepackage{hyperref}
\usepackage{booktabs}
\usepackage{tabu}
\usepackage{longtable}
\usepackage{lscape} 
\usepackage{rotating}

\newtheorem{example}{Example}

\allowdisplaybreaks

\makeatletter
\def\ps@pprintTitle{%
  \let\@oddhead\@empty
  \let\@evenhead\@empty
  \def\@oddfoot{\reset@font\hfil\thepage\hfil}
  \let\@evenfoot\@oddfoot
}
\makeatother

\begin{document}

\begin{frontmatter}
\title{Hub location with congestion and time-sensitive demand}

\author[uca]{Carmen-Ana Domínguez-Bravo\corref{cor1}}
\ead{carmenana.dominguez@uca.es}

\author[uca]{Elena Fern\'{a}ndez}
\ead{elena.fernandez@uca.es}

\author[unab]{Armin L\"{u}er-Villagra}
\ead{armin.luer@unab.cl}

\address[uca]{Statistics and Operations Research Department, Universidad de Cádiz, Puerto Real, Spain}

\address[unab]{Department of Engineering Sciences, Universidad Andres Bello, Santiago, Chile}

\cortext[cor1]{Corresponding author}

\date{}

\begin{abstract}
This work studies the effect of hub congestion and time-sensitive demand on a hub-and-spoke location/allocation system.
The Hub Location with Congestion and Time-sensitive Demand Problem is introduced, which combines these two main characteristics. On the one hand, hubs can be activated at several service levels, each of them characterized by a maximum capacity, expressed as the amount of flow that may circulate through the hub, which is associated with a hub transit time. On the other hand, alternative levels are available for served commodities, where each demand level is characterized by its amount of demand, unit revenue, and maximum service time. In this problem the efficiency of a hub-and-spoke system is given by the maximum net profit it may produce.
To the best of our knowledge this is the first work where hub congestion and time-sensitive demand are jointly considered.\\
Two alternative mixed-integer linear programming formulations are proposed. They include a new set of constraints, which are necessary to guarantee the consistency of the obtained solutions under the presence of the \textit{capacity-type} constraints derived from hub service levels and served demand levels. The efficiency of the formulations is analyzed through a set of computational experiments. The results of the computational experiments allow to study the structure of the obtained solutions and to analyze how the different parameters affect them.

\end{abstract}

\begin{keyword}
Location \sep%
Hub Location \sep
\end{keyword}

\end{frontmatter}


\section{Introduction}

Hub location problems (HLPs) aim to locate a special kind of facilities called hubs, used for consolidation, sorting, and classification of flows connecting given origin/destination (OD) demand pairs, commonly named commodities.

\textit{Hub-and-spoke} are among the most popular topologies for networks built around hubs, where flows are routed through paths that traverse at least one hub. The design and operation of hub-and-spoke networks is motivated by the requirement of a smaller number of links and the opportunity of achieving economies of scale by consolidating flows between hubs. Hub-and-spoke networks are commonly designed using hub location optimization models.

The first `modern' problem statement and mathematical formulations related to hub location models were presented by \citet{OKelly1986,OKelly1987}, whereas  \citet{Campbell1992,Campbell1994} introduced the `hub equivalent' statement and linear formulations of the four main classes of location models, i.e. covering, median, center, also for the cases with fixed costs. Thereafter, hub location flourished as a research area, due to two main reasons. On the one hand, because of the importance for the transportation industries, such as cargo, air passenger transportation, public transport, liner shipping, and telecommunications. On the other hand, hub location has attracted many researchers because of its theoretical interest and the wide range of methodological alternatives for addressing this type of problems.

 Mathematical formulations modeling hub-and-spoke systems are usually challenging to solve as they involve joint location, allocation, and routing decisions. Typically, they require specialized solution methods, such as Benders decomposition \citep{Gelareh2008proc, DeCamargo2009, Contreras2011a, DeCamargo2011, DeSa2013, OKelly2015, Merakli2017, Ghaffarinasab2019, Taherkhani2020, Taherkhani2021, Ghaffarinasab2022}, Lagrangean methods \citep{Elhedhli2005, Elhedhli2010, Tanash2017, Alibeyg2018, Ghaffarinasab2022}, metaheuristics \citep{Calik2009, Mohammadi2011, LV2013, Silva2017, Azizi2018, LV2019}, branch\&cut \citep{Rothenbacher2016, Zetina2017, Ghaffarinasab2019}, to name a few.

 Initial \textit{fundamental} models have evolved to multiple extensions capturing more realistic aspects, such as competitive models \citep{Marianov1999, LV2013,Sasaki2014}, more accurate cost structures \citep{Bryan1998, OKelly1998, Kimms2006, LV2019}, multiperiod problems \citep{Gelareh2015, Alumur2016, Correia2018, Monemi2021}, hub congestion \citep{Marianov2003, Mohammadi2011, Ishfaq2012, Rahimi2016, Karimi2020}, accurate demand modeling \citep{LV2013, ArosVera2013, OKelly2015, Alibeyg2016, Taherkhani2020}, or profit maximizing problems \citep{Alibeyg2016, Alibeyg2018, Taherkhani2019, Taherkhani2020, Taherkhani2021}. Frequent literature reviews and book chapters show the permanent interest within the Location Analysis community \citep{Campbell2002, Alumur2008, Kara2011, Campbell2012, Farahani2013,  Contreras2015, Alumur2021}.

Among the topics that have received increasing attention in the last years we can mention congestion and time-sensitive demand. Timely service is crucial in transportation networks so congestion due to peaks of demand should be avoided because of its multiple obnoxious effects. The most obvious consequence of congestion is the increase in expected travel times, leading to delays, missed connections, etc. In its turn, these delays may discourage users, producing a negative impact on demand and, eventually, on the overall revenue attained by the transportation network. Hence these two topics are closely related although, to the best of our knowledge, they have not yet been addressed jointly.  This is precisely the focus and main contribution of this paper, where we introduce and solve a profit maximization hub network design problem, where hub congestion and time-sensitive demand are jointly considered.

Profit maximization HLPs assume that existing demand needs not be necessarily served and that served commodities produce a revenue. In these problems, the decision-maker has to extend the set of decisions so as to also determine the commodities that must be served. The objective is now to maximize the total net profit, i.e., the difference between the total revenue for the served commodities minus the usual (activation + routing) costs. The interested reader is addressed to  \citep{Alibeyg2016,Alibeyg2018,Taherkhani2019,Taherkhani2020,Taherkhani2021}.

We model congestion by relating transit times at hubs to the amount of flow circulating through each of them. Furthermore, we model time-sensitive demand by relating the service demand of each commodity to its overall travel time through the activated hub-and-spoke network as well as to the revenue it produces to the owner (manager) of the network if it is served. For this, the overall travel times of service paths depend, not only on the traversal times of their arcs, but also on the capacity-dependent transit times at the hubs that are crossed. These overall travel times affect both unit revenues and served demand. Specifically, we assume that the revenue produced by each served commodity decreases as its overall travel time increases. Furthermore, we assume that demand is also sensitive to service route attributes. In particular, we assume that the amount of demand of each commodity is a (simple) concave function of its travel time. Thus, \textit{most} of the demand will correspond to \textit{standard} users who accept \textit{intermediate} service times and produce \textit{intermediate} revenues, whereas some (fewer) users will impose \textit{express} (small) service times but would produce higher revenues and some other users will admit longer travel times associated with smaller revenues to the company.

  We are not aware of any existing work regarding a hub network design problem considering hub congestion together with time-sensitive demand. In fact, to our knowledge, the modeling approach that we follow in which demand is sensitive to routes duration has not been previously adopted in the context of profit maximizing HLPs. In this paper we introduce a new problem that we call the Hub Location with Congestion and Time-sensitive Demand Problem (HLCTDP), where these two, closely related, issues are jointly addressed. We study some properties of this novel and challenging hub network design problem, and develop two mathematical programming formulations, which are empirically compared through a series of computational experiments. The most efficient formulation is further enhanced by exploiting some optimality conditions. This allows us to solve to proven optimality instances with up to 50 nodes in reasonable computing times.

The remainder of the paper is organized as follows. Section \ref{literature} reviews the most relevant related literature. Section \ref{problem} gives the formal definition of the HLCTDP, and Section \ref{mathematicalformulation} presents two mathematical formulations that allow solving it with off-the-shelf solvers. In Section \ref{optimalityconditions}, noting the difficulty to solve the formulations, we study some properties of the problem, leading to some optimality conditions that allow reducing the number of decision variables.
 Section \ref{experimentsresults} explains how we have generated the HLCTDP test instances used, describes the computational experiments, and presents and analyzes the obtained results. The paper ends in Section \ref{conclusion} with some conclusions and comments of possible avenues for future research.

 \section{Literature review}\label{literature}
In this section we review the most relevant literature on the topics related to this paper.

 \subsection{HLPs with congestion}

 Modeling congestion in HLPs has generated permanent interest within the hub location literature. The approaches used in this context are multiple and diverse, but can be grouped in: queueing models, models with a non-linear term representing congestion, and other models. Table \ref{table:congestion} compares relevant literature with our work regarding congestion in HLPs.

 To the best of our knowledge, the first HLP considering congestion was introduced by \cite{Marianov2003}, where the authors model hubs as M/D/c queueing systems, leading to a family of chance capacity constraints, and solve the model using a tabu search-based heuristic.
 After that, \cite{Elhedhli2010} extended previous work by adding capacity decisions to the problem. They model the hub network as M/M/1 queues and solve the resulting model using Lagrangean methods, in a similar way to \citep{Elhedhli2005}.

 Later, \cite{Mohammadi2011} used an $M/M/c$ queue model to represent congestion for a hub covering location problem. To capture the effect of limited capacity at hubs and its impact on network design, \citet{Ishfaq2012} model hub operations as a $GI/G/1$ system.

  More recently, \cite{Rahimi2016} proposed a multi-objective model, where hub congestion is modeled using $M/M/c/k$ systems. Closely related to the work of \cite{Ishfaq2012}, \cite{Karimi2020} use a $GI/G/c$ queueing model to capture the congestion at hubs and links, and solve the resulting model using a hybrid metaheuristic.

 From a different perspective, a body of literature adds a penalty term to the objective function to take congestion into account. The first model following this approach was proposed in \citep{Elhedhli2005}, where congestion is included as a non-linear cost term. The authors then linearize the resulting model, and use Lagrangean methods to obtain feasible solutions. Similarly, \cite{DeCamargo2009} represent congestion as a convex cost function and use a generalized Benders decomposition method to solve instances to optimality. They also establish the trade-off between the cost saving due to economies of scale and the congestion generated. Later, \cite{DeCamargo2011} proposed a hybrid approach, combining outer-approximation and generalized Benders cuts to solve the single-allocation variant of the problem. More recently, \cite{Kian2016} used power-law congestion costs, leading to a conic quadratic model. The authors compare two ways of modeling and solving the problem, highlighting pros and cons.

  Finally, other works propose different approaches to model congestion. In \citep{Yang2015} the authors use a two-stage stochastic programming model with recourse to consider stochastic demand and hub congestion. Multiple discrete levels of hub capacity and congestion are considered in \citep{Alumur2018}, leading to linear models, which can be solved by off-the-shelf solvers. Here we also consider a discrete set of \textit{congestion levels}, which are determined by some lower and upper limits on the overall flow at the hubs.

\begin{table}[htbp]\footnotesize
\caption{HLP with congestion.}
\label{table:congestion}
\centering
\begin{tabular}{lcccccccc}
\toprule
\textbf{Authors} & \textbf{Approach} & \textbf{Method} & \textbf{H. Cap.} & \textbf{Demand} & \textbf{Congestion} & \textbf{MILP} & \textbf{OF} & \textbf{Alloc}. \\
\midrule
\cite{Elhedhli2005}  & A    & Lag  & --   & 1L     & Convex  & -- & mC       & M \\
\cite{DeCamargo2009} & E    & BD   & 1L   & 1L     & Convex  & -- & mC       & M \\
\cite{Elhedhli2010}  & A    & Lag  & --   & 1L     & M/M/1   & -- & mC       & S \\
\cite{DeCamargo2011} & E    & BD   & --   & 1L     & Convex  & -- & mC       & S \\
\cite{Mohammadi2011} & A    & MaH  & Prob & 1L     & M/M/c   & \checkmark & mC       & S \\
\cite{Ishfaq2012}    & E    & MH   & --   & MaxST  & GI/G/1  & -- & mC       & M \\
\cite{Rahimi2016}    & A    & MH   & ML   & 1L     & M/M/c/K & -- & mC/mMT   & S \\
\cite{Alumur2018}    & E    & Form & ML   & MaxST  & ML      & \checkmark & mC       & S, M \\
\cite{Azizi2018}     & E, A & MH   & ML   & Stoch  & M/G/1   & \checkmark & mC       & S \\
\cite{Karimi2020}    & A    & MaH  & --   & 1L     & GI/G/1  & -- & mC/mMT  & M \\
\midrule
This work (2023) & E & Form & ML & Time, ML & ML & \checkmark & MP & M \\
\bottomrule
\end{tabular}

\begin{flushleft}
\underline{\textbf{Notation:}}\\
\textbf{Approach:} A: approximate; E: Exact; \\
\textbf{Method:} BD: Benders Decomposition; Form: Formulation; Lag: Lagrangean; MaH: Math-heuristic; MH: Metaheuristic;\\
\textbf{H. Cap. (Hub Capacity)}: 1L: One level; ML: Multiple levels; Prob: Probabilistic;\\
\textbf{Demand (Demand Modeling)}: MaxST: Maximum Service Time; Stoch.: Stochastic; Time: Time Sensitive demand;\\
\textbf{Congestion (Hub Congestion Modeling)}: Convex: Convex Function; ML: Multiple levels;\\
\textbf{MILP (Mixed Integer Linear Programming):} \checkmark: Yes; --: No;\\
\textbf{OF (Objective Function)}: mC: Minimize Total Cost; mMT: Minimize Maximum Time; MP: Maximize Profit;\\
\textbf{Alloc. (Allocation)}: S: Single Allocation; M: Multiple Allocation.
\end{flushleft}
\end{table}

\subsection{HLPs with sensitive demand}

 Modeling demand in HLPs is one of the most popular topics within the hub location community. The most studied perspectives are demand sensitivity and demand uncertainty. Table \ref{table:sensitive_demands} compares previous work with our proposal  regarding HLPs with sensitive demand.

Most classical HLPs consider fixed demand, i.e. they assume it is deterministic and it does not depend on other factors. The only exception among deterministic models being hub covering and center models \citep{Campbell1994}, where users demand is captured if their service routes do not exceed a maximum time or cost threshold \citep[see, e.g.,][]{Kara2003, Tan2007, Wagner2008, Calik2009, Peker2015, Silva2017}.

Further extensions aiming to model sensitive demand can be classified in proportional choice models, like Huff/gravitational and logit models, and deterministic choice models. To our knowledge, \cite{Eiselt2009} present the first HLP where demand sensitivity is expressed through a discrete choice model. The authors study the problem faced by an entrant company aiming to maximize its market capture, assuming that customers demand can be modeled with a gravity-like choice model, solving it using heuristic concentration.

Some examples of works studying hub-and-spoke systems for different application areas where logit models have been used to model customers' choice are \citep{LV2013} for airlines, \citep{ArosVera2013} for park-and-ride facilities, and \citep{Mahmoodjanloo2020} for postal services.

 \cite{OKelly2015} address a HLP with price-sensitive demand. The model maximizes the social net benefit, and there are three parallel types of services, each of them with different associated demand: direct, 1-hub-stop and 2-hub-stops. Two formulations are proposed, as well as a Benders decomposition solution algorithm, together with a detailed analysis of the price equilibrium associated with the solutions obtained.

 From a different perspective, \cite{Alibeyg2016} consider some profit maximizing HLPs with multiple demand levels, where each level implies a certain demand and revenue, at most one of the demand levels can be selected simultaneously, and hubs do not have capacity constraints.

 \cite{Taherkhani2020} also address a profit maximizing HLP with multiple demand classes, which incorporates capacity constraints on the hubs. They provide deterministic and stochastic versions of the problem, which are solved with Benders decomposition and Sample Average solution algorithms in the deterministic and stochastic versions, respectively.

\cite{Alibeyg2016} and \cite{Taherkhani2020} associate each demand level with a revenue. However, neither of these works assume that demand is sensitive to service route attributes. On the contrary, in our work we assume that the users demand is sensitive to the characteristics of their service routes, namely to the overall duration of the routes. Moreover, since we also consider hub congestion, and the duration of routes is affected by it, the users' demand is sensitive to the congestion of the hub-and-spoke system.

  Multiple approaches have been followed in HLPs considering demand uncertainty. Some examples are stochastic programming \citep{Contreras2011b, Yang2015}, stochastic processes \citep{Azizi2018}, and robust optimization \citep[see, e.g.,][]{Merakli2016, Merakli2017, Zetina2017, Wang2020}. In this work we do not assume uncertainty, as we want to keep the problem solvable, as it already incorporates other complexities, including profit maximization, hub congestion, and demand sensitivity.

 \begin{table}[htbp]\footnotesize
\caption{HLPs with sensitive demands.}
\label{table:sensitive_demands}
\begin{tabular}{lcccccc}
\toprule
\textbf{Authors} & \textbf{Method} & \textbf{H Cap.} & \textbf{Demand} & \textbf{Congestion} & \textbf{OF} & \textbf{Alloc.} \\
\midrule
\cite{Eiselt2009}        & MH            & -- & AF     & -- & MD & M \\
\cite{Contreras2011b}    & BD, MH         & -- & Stoch & -- & mC & M \\
\cite{Alumur2012}        & Form, SP, MH & -- & UD     & -- & mC & S,M \\
\cite{Ishfaq2012}        & Lagr, MH     & -- & MaxST  & GI/G/1 & mC & M \\
\cite{LV2013}            & MH            & -- & Logit  & -- & MP & M \\
\cite{OKelly2015}        & BD            & -- & Price  & -- & MP & M \\
\cite{Peker2015}         & Form         & -- & Cost   & -- & MD & S, M \\
\cite{Yang2015}          & SP            & 1L & Stoch & -- & mC & M \\
\cite{Alibeyg2016}       & Form         & ML &  ML  & -- & MP & M \\
\cite{Merakli2016}       & BD, RO        & -- & UD     & -- & mMT & M \\
\cite{Merakli2017}       & BD, RO        & 1L & UD     & -- & mMT & M \\
\cite{Zetina2017}        & Form, BC     & -- & UD     & -- & mC & M \\
\cite{Correia2018}       & Form         & 1L & Stoch & -- & mC & M \\
\cite{Alumur2018}        & Form         & ML & MaxST  & ML & mC & S, M \\
\cite{Azizi2018}         & MH            & ML & Stoch & M/G/1 & mC & S \\
\cite{Mahmoodjanloo2020} & MaH           & -- & Logit & -- & MP & S \\
\cite{Taherkhani2020}    & BD, SP        & 1L & ML, Stoch & -- & MP & M \\
\cite{Wang2020}          & RO            & 1L & UD     & -- & mC & M \\
\cite{Monemi2021}        & BD            & -- & Serial  & -- & mC & S \\
\cite{Taherkhani2021}    & BD, SP        & 1L & Revenue & -- & MP & M \\
\cite{Ghaffarinasab2022} & BD, Lag       & -- & Bernoulli & -- & mC & S, M \\
\midrule
This work (2023)        & Form          & ML & Time, ML & ML & MP & M \\
\bottomrule
\end{tabular}
\begin{flushleft}
\underline{\textbf{Notation:}}\\
\textbf{Method:}
BC: Branch\&Cut;
BD: Benders Decomposition;
Form: Formulation;
Lag: Lagrangean;
MaH: Math-heuristic;
MH: Metaheuristic;
RO: Robust Optimization;
SP: Stochastic Programming;
\\
\textbf{H Cap. (Hub Capacity):} 1L: Single Level; ML: Multiple Levels;\\
\textbf{Demand (Demand Modeling):}
AF: Attraction Functions;
Bernoulli: Bernoulli random;
Cost: Cost sensitive;
Logit: Logit Discrete Choice Model;
MaxST: Maximum Service Time;
Price: Price Sensitive;
Revenue: Revenue sensitive;
Serial: Serial demand;
Stoch: Stochastic;
UD: Uncertain Demand;
\\
\textbf{Congestion (Hub Congestion Modeling)}:
ML: Multiple levels;
A/B/C: Conventional queueing notation;\\
\textbf{OF (Objective Function):}
mC: Minimize Total Costs;
MD: Maximize Demand Captured;
MP: Maximize Profit;
mMT: Minimize Maximum Service Time;\\
\textbf{Alloc. (Allocation)}:
S: Single Allocation;
M: Multiple Allocation.
\end{flushleft}
\end{table}

\subsection{Profit maximizing HLP}
Profit maximizing HLPs incorporate decisions on which commodities to serve to those related to the design of the hub network.  Unlike HLPs where all demand must be served and decisions are based on the time or cost of the offered routes, in profit maximizing HLPs the decisions on the commodities to serve rely on the net profit that service routes would provide to the decision maker. Table \ref{table:profit_maximizing} compares relevant literature on this topic with our work.

 To our knowledge, \citep{Alibeyg2016} is the first work regarding a profit maximizing HLP. They introduce as base model a specific hub network design problem with profits, which is later extended by considering service commitments and direct connections. They further extend the model to include setup costs in access or bridge edges, and multiple demand levels.

 Later, \cite{Alibeyg2018} proposed exact algorithms for two profit maximizing HLPs based on Lagrangean methods, together with reduction tests and partial enumerations. \cite{Taherkhani2019} also studied profit maximizing HLPs, focusing their study on the allocation strategies, i.e., how non-hub nodes are allocated to hub nodes, and on modeling direct connections between non-hub nodes. The resulting models are solved using CPLEX after performing variable fixing tests.

 In \citep{Taherkhani2020} a Benders decomposition solution algorithm is proposed for a capacitated profit maximizing HLP with multiple demand classes. After solving a deterministic model, the authors extend their procedures to cope with a stochastic programming formulation. The numerous computational improvements developed allow the authors to solve to proven optimality instances with up to 500 nodes and three demands classes.
 Extending their own work, in \citep{Taherkhani2021} the authors introduce some robust stochastic profit maximizing HLPs, considering both stochastic demand and uncertain revenues. They develop different algorithmic improvements allowing to solve large instances.

 In this work, we state and solve a profit maximizing HLP. Instead of modeling uncertainty, we focus on the dependence between time-sensitive demand and hub congestion.

\begin{table}[htbp]\footnotesize
\caption{Profit maximizing HLPs.}
\label{table:profit_maximizing}
\begin{tabular}{lccccccc}
\toprule
\textbf{Authors} &
\textbf{Method} &
\textbf{H Cap.} &
\textbf{Demand} &
\textbf{Congest.}   &
\textbf{MILP} &
\textbf{Comp.} &
\textbf{Alloc.}\\
\midrule
\cite{LV2013}            & MH     & -- & Logit     & -- & -- & \checkmark & M      \\
\cite{OKelly2015}        & BD     & -- & Price     & -- & \checkmark & -- & M      \\
\cite{Alibeyg2016}       & Form   & ML & ML        & -- & \checkmark & -- & M      \\
\cite{Alibeyg2018}       & Lagr   & -- & 1L        & -- & \checkmark & -- & M      \\
\cite{Taherkhani2019}    & Form   & -- & 1L        & -- & \checkmark & -- & S,M,rA \\
\cite{Mahmoodjanloo2020} & MH     & -- & Logit, ML & -- & -- & \checkmark & S      \\
\cite{Taherkhani2020}    & BD, SP & 1L & ML        & -- & \checkmark & -- & M      \\
\cite{Taherkhani2021}    & BD, SP & 1L & Revenue   & -- & \checkmark & -- & M      \\
\midrule
This work (2023)        & Form & ML & Time, ML & ML & \checkmark & -- & M \\
\bottomrule
\end{tabular}
\begin{flushleft}
\underline{\textbf{Notation:}}\\
\textbf{Method:} BD: Benders Decomposition; Form: Formulation; Lag: Lagrangean; MH: Meta-Heuristic; SP: Stochastic Programming;\\
\textbf{H Cap. (Hub Capacity):} 1L: Single Level; ML: Multiple Levels;\\
\textbf{Demand (Demand Modeling)}: 1L: Single Level; ML: Multiple Levels; Logit: Logit discrete choice model, Price: Price sensitive demand; Time: Time sensitive demand;\\
\textbf{Congest. (Hub Congestion modeling)}: ML: Multiple Levels; --: No; \\
\textbf{MILP (Mixed Integer Linear Programming):} \checkmark: Yes; --: No; \\
\textbf{Comp. (Competitive):} \checkmark: Yes; --: No;\\
\textbf{Alloc. (Allocation):} S: Single Allocation; M: Multiple Allocation; rA: $r$-Allocation.
\end{flushleft}
\end{table}

 \section{Definition of the Hub Location with Congestion and Time-sensitive Demand Problem}
 \label{problem}

\begin{figure}[htbp]
\begin{subfigure}{0.50\linewidth}
	\begin{center}
		\includegraphics[scale=0.5]{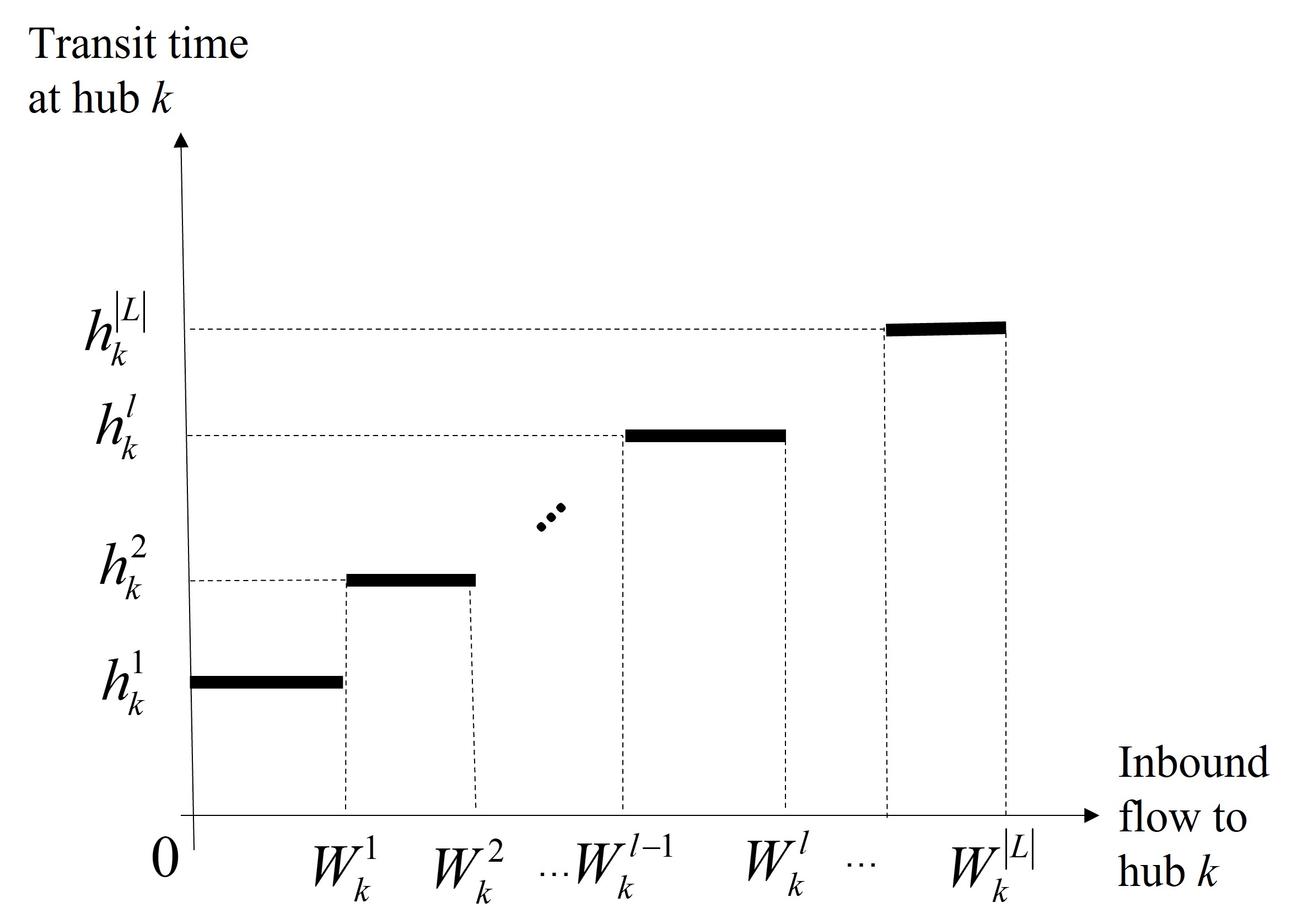}
	\end{center}
	\caption{Inbound flows and transit times at potential hubs.}\label{Fig:HubTransitTime}
\end{subfigure}
\begin{subfigure}{0.50\linewidth}
\includegraphics[scale=0.5]{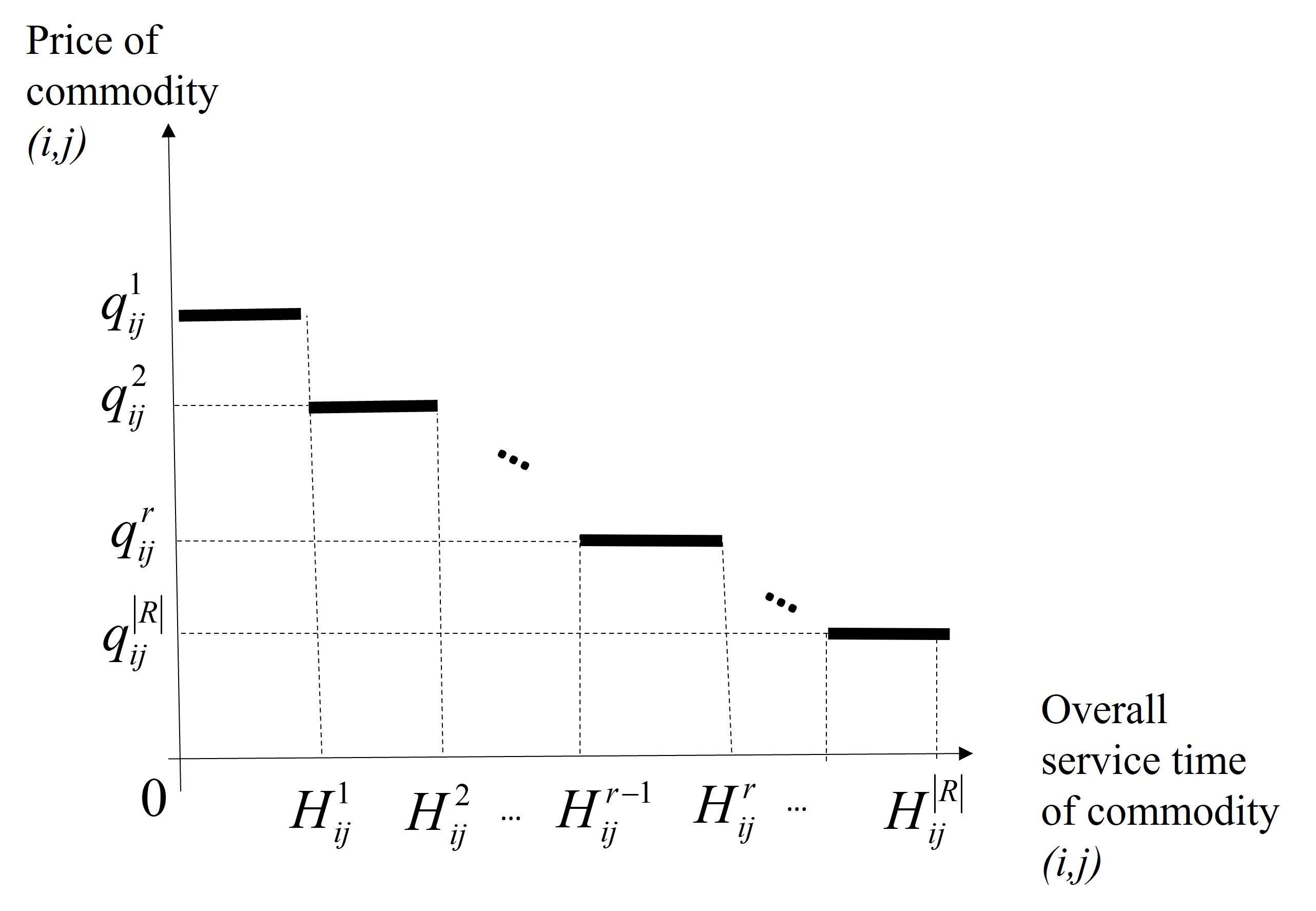}
\caption{Overall service time and price charged to the users, for a given commodity.}\label{Fig:price}
\end{subfigure}
\caption{Relationship between selected parameters.}
\end{figure}

We consider a complete directed graph $G = (V, A)$, where $V=\{1, 2,\dots, n\}$ represents the set of users, and $K\subseteq V$ the set of potential locations for the hubs. For each $(i,j)\in A$, $t_{ij}$ and $c_{ij}$ respectively denote the unit travel time and the unit routing (transportation) cost through arc $\left(i,j\right)\in A$.
Let $L$ denote the index set of possible service levels at potential hubs. Specifically, associated with each potential hub $k\in K$ and service level $l\in L$, there are two quantities $h_k^l$ and $W_k^l$, which respectively denote the transit time and the maximum inbound flow at hub $k$, if it operates at service level $l$. Figure~\ref{Fig:HubTransitTime} illustrates the relationship between transit times ($h_k^l$) and maximum inbound flows ($W_k^l$) at a given hub node $k$. As can be seen, we model this relationship as an increasing stepwise function. That is, we assume that service level $l$ corresponds to an inbound flow in the interval $(W_k^{l-1}, W_k^l]$ where $W_k^0=0$. Higher inbound flows imply higher transit times, although we assume a constant transit time for each service level.

Let $\mathcal{C}\subset V\times V$ denote a given set of commodities, i.e. ordered pairs of nodes with service demand. Let also $R$ denote the index set for the commodities demand levels. Associated with each commodity $(i,j)\in \mathcal{C}$, and demand level $r\in R$ we are given the following quantities that apply if commodity $(i, j)$ is served at demand level $r$: ($i$) $w_{ij}^r$, which denotes its amount of demand, i.e. the flow that must be sent from $i$ to $j$;
($ii$) $H_{ij}^{r}$, an upper bound on the overall service time (travel plus transit at the traversed hubs) from $i$ to $j$, and, ($iii$) $q_{ij}^{r}$, the prize charged to the user (revenue obtained by the decision maker) per unit of routed flow.
Figure~\ref{Fig:price} illustrates the relationship between the price charged to the user ($q_{ij}^{r}$) and the maximum service time ($H_{ij}^{r}$), for a given commodity $(i,j)\in \mathcal{C}$. We model this relationship as a decreasing stepwise function, so a larger service time implies a cheaper user fee, although the fee remains constant for a service time within the interval $(H_{ij}^{r-1}, H_{ij}^{r}]$, $r\in R$, where $H_{ij}^{0}=0$.


Commodity $(i, j)\in \mathcal {C}$ may not be served at all. In that case it produces no revenue nor incurs any cost. When commodity $(i, j)$ is served it will be associated with a unique demand level $r\in R$.  This will produce a revenue of $q_{ij}^r w_{ij}^{r}$ and will incur a routing cost, which depends both on the amount of demand $w_{ij}^r$ and on the path used for routing it.

The routing path for the flow $w_{ij}^r$ must be of the form $i-k-m-j$
where $k, m\in K$ must be open hubs and it is possible that $k=m$. In that case the routing path will be of the form $i-k-j$.
When either $i$ or $j$ are open hubs, routing paths reduce to connections of the form $i-m-j$ or $i-k-j$.
When both $i$ and $j$ are open hubs, the routing path is the direct connection $i-j$.
The routing cost for serving commodity $(i, j)$ at demand level $r\in R$ through path $i-k-m-j$, is $w_{ij}^r(c_{ik}+\alpha\,c_{km}+c_{mj})$, where $\alpha$ is the interhub routing cost discount factor.

The condition that the service time of the commodity cannot exceed $H_{ij}^r$ if it is served at demand  level $r\in R$, is modeled by considering the possible cases.
When $k\ne m$ and hubs $k$ and $m$ are activated at service levels $l, l'\in L$, respectively, the following inequality must be satisfied: $t_{ik}+h_k^l+\gamma\,t_{km}+h_{m}^{l'}+t_{mj}\le H_{ij}^r$, where $\gamma$ is the interhub routing time discount factor.  Note that this constraint reduces to $h_i^l+\gamma\,t_{im}+h_{m}^{l'}+t_{mj}\le H_{ij}^r$ and $t_{ik}+h_{k}^l+\gamma\,t_{kj}+h_j^l\le H_{ij}^r$, respectively, when $i$ or $j$ are open hubs. When both $i$ and $j$ are open hubs, the service time  is $h_i^l+\gamma\,t_{ij}+h_{j}^{l'}$. When $k= m$ and hub $k$ is activated at service level $l\in L$,  the time limit inequality becomes $t_{ik}+h_k^l+t_{kj}\le H_{ij}^r$.

As explained, if hub $k\in K$ is activated at service level $l$, the total demand of the commodities routed through $k$ cannot exceed the upper bound $W_k^l$. Activating hub $k$ at service level $l$ incurs a fixed cost $G_{k}^{l}$.

In the HLCTDP the following decisions must be made: $(i)$ Select the set of hubs to open, together with their associated service levels; ($ii$) Select the commodities to serve, and their associated demand levels; and, $(iii)$ Determine for each of the served commodities a routing path through at most two intermediate hub nodes, respecting the congestion limitations at the hub nodes and the maximum service times of the demand levels of the served commodities.
The objective is to maximize the difference of the total profit for the served commodities minus the sum of the setup costs of the activated hubs and the routing costs of the served commodities.

The HLCTDP has as particular cases several well-known hub location models. When $|L|=1$, so there is one single service level at the hubs, the HLCTDP is a particular case of a hub network design problem with profits \citep{Alibeyg2016}. When $|R|=1$, so there is one single service demand level for the commodities, the HLCTDP becomes a hub location model with congestion \citep[see, e.g.][]{Elhedhli2005,Elhedhli2010,Kian2016,Alumur2018,Luer-Villagra2019}. When both $|R|=1$ and $|L|=1$ and profits are sufficiently large, the HLCTDP reduces to the classical uncapacitated $p$-hub median problem \cite{Campbell1994}.  Therefore, the HLCTDP is clearly NP-hard.

Throughout we use the notation:
\begin{itemize}
	\item $C_{ijkm}=c_{ik}+\alpha\,c_{km} +c_{mj}$: unit routing cost for sending commodity $(i,j)\in\mathcal{C}$ via  $k$ and $m$.
    \item $\overline C_{ijkm}^{r}=w_{ij}^r\left(q_{ij}^r-(c_{ik}+\alpha\, c_{km}+c_{mj})\right)$: net profit for serving commodity $(i,j)\in\mathcal{C}$ at service level $r\in R$ and routing it via hubs $k$ and $m$.
\item $T_{ijkm}=t_{ik}+\gamma\,t_{km}+t_{mj}$: the travel time from $i$ to $j$ via hubs $k$ and $m$, without taking into account transit times at the intermediate hubs.
\end{itemize}

\section{Mathematical programming formulations}\label{mathematicalformulation}
In this section we introduce two mixed-integer linear programming (MILP) formulations for the HLCTDP. The main target is to develop formulations that allow solving instances of reasonable sizes to proven optimality. For this it is necessary that the formulations produce sufficiently tight bounds so that the size of the enumeration tree does not \textit{explode} while not being too demanding in terms of memory requirements. This is a challenge since the HLCTDP needs to incorporate additional decisions to a \textit{classical} hub location problem, whose most effective formulations already use 4-index routing variables for determining the routing paths (two indices for the origin/destination of the involved commodity plus two additional ones for the intermediate hubs in the routing path).

The two formulations that we introduce use the same sets of decision variables to model the strategic decisions of the HLCTDP, namely the location of hubs to activate together with their service levels, as well as the choice of the commodities to serve together with their demand levels.
For the location decisions we adopt the natural option by incorporating one additional index to usual location variables to indicate the service level of the hubs. Specifically, the location decision variables that we use are:

\begin{itemize}
	\item $Y_{k}^{l}$: binary variable that takes the value 1 if and only if a hub is opened at node $k \in K$ at service level $l\in L$, and 0 otherwise. There are $\mathcal{O}(|K|\times|L|)$ variables $Y$.
\end{itemize}

For determining the served commodities and associated demand levels we use the following variables:
\begin{itemize}
 \item $\beta_{ij}^{r}$: binary variable that takes the value 1, if commodity $\left(i,j\right) \in \cal C$ is served at demand level $r\in R$; 0, otherwise. There are $\mathcal{O}(|\mathcal C|\times|R|)$ variables $\beta$.
\end{itemize}

The essential difference among the two formulations is how the variables that determine the routing paths, and their associated flows, are defined.
Note that, contrary to other HLPs where the routing variables can be relaxed to take continuous values, in the HLCTDP it is necessary to impose explicitly that routing variables take binary values. This happens because the maximum service times (travel plus transit) of the served commodities and the maximum flow limitation at the activated hubs lead to knapsack-type constraints. This is yet another reason that makes the number of routing variables particularly important here.
Therefore, in order to limit the overall number of routing decision variables, we exclude the most intuitive possibility of incorporating three additional indices to classical routing variables: one to indicate the demand level of the involved commodity, and (up to) two more to indicate the service level of the traversed hubs.

The first formulation that we introduce uses \textit{traditional} 4-index routing variables, at the expenses of needing additional sets of decision variables to determine the flows circulating at the arcs so as to control the congestion at the activated hubs. The second formulation uses 5-index routing variables, which allow to easily identify the amount of flow circulating through each arc.

\subsection{Formulation F1: 4-index routing variables}

 In this formulation the routing, and additional decision variables are the following:
 \begin{itemize}
  \item $x_{ijkm}$: 1, if commodity $\left(i,j\right) \in \cal C$ is routed through hubs $k,m \in K$ (in that order); 0, otherwise.
  \item $f_{ijk}^{l}$: 1, if hub $k \in K$ at service level $l\in L$, is the first hub in the OD path of commodity $\left(i,j\right) \in \cal C$; 0, otherwise.
  \item $s_{ijk}^{l}$: 1, if hub $k \in K$ at service level $l\in L$,  is the second hub in the OD path of commodity $\left(i,j\right) \in \cal C$; 0, otherwise.
  \item $o_{ijk}^{l}$: 1, if hub $k \in K$ at service level $l\in L$,  is the only hub in the OD path of commodity $\left(i,j\right) \in \cal C$; 0, otherwise.
   \item $g_{ijkm}$: flow of commodity $\left(i,j\right) \in \cal C$ in an OD path through hubs $k,m\in K$.
 \end{itemize}

There are $\mathcal{O}\left(|\mathcal{C}|\times|K|^2\right)$  variables $x$, $\mathcal{O}\left(|\mathcal{C}|\times|K|\times|L|\right)$  variables of each of the types $f$, $s$, and $o$, and a number of $g$ variables $\mathcal{O}\left(|\mathcal{C}|\times|K|^2\right)$. Formulation F1 is:

\begin{small}
\begin{align}
F1\quad  \max &\sum_{\left(i,j\right) \in \cal C} \sum_{r \in R} q_{ij}^{r} w_{ij}^{r} \beta_{ij}^{r} - \sum_{\substack{k \in K, l \in L}} G_{k}^{l} Y_{k}^{l} &&  \label{of}\\
\text{s.t.}&  \sum_{l\in L} Y_{k}^l \le 1 \quad && k\in K   \label{1oneFlowLevel}\\
&  \sum_{r\in R} \beta_{ij}^{r} \leq 1 \qquad && \left(i,j\right)\in \cal C   \label{demandlevels}\\
&  \sum_{r \in R} \beta_{ij}^{r}= \sum_{k,m \in K} x_{ijkm} \qquad &&  \left(i,j\right)\in \cal C  \label{demandLevelPath}\\
&  \sum_{l \in L} f_{ijk}^{l} = \sum_{m \in K: m \neq k} x_{ijkm} \qquad &&  \left(i,j\right)\in {\cal C}, k \in K   \label{yandxBound}\\
&  \sum_{l \in L} s_{ijk}^{l} = \sum_{m \in K: m \neq k} x_{ijmk} \qquad &&  \left(i,j\right)\in {\cal C}, k \in K  \label{yoverlineandxBound}\\
&  \sum_{l \in L} o_{ijk}^{l} = x_{ijkk}  \qquad &&  \left(i,j\right)\in {\cal C}, k \in K   \label{yhatandxBound}\\
& 	x_{ijkk} + \sum_{m\in K: m \neq k} \left(x_{ijkm} + x_{ijmk}\right) \leq \sum_{l \in L} Y_{k}^{l}  \qquad &&  \left(i,j\right)\in {\cal C}, k \in K
  \label{xAndzBound}\\
&  f_{ijk}^{l} + s_{ijk}^{l} + o_{ijk}^{l} \leq Y_{k}^{l} \qquad &&  \left(i,j\right)\in {\cal C}, k \in K, l \in L \label{yAndYBound}\\
&  g_{ijkm} \leq \left(\max_{r \in R} w_{ij}^r \right)\,x_{ijkm} \qquad &&  \left(i,j\right)\in {\cal C}, k,m \in K  \label{gxbound} \\
&  \sum_{k,m \in K} g_{ijkm} = \sum_{r \in R} w_{ij}^{r} \beta_{ij}^{r} \qquad && \left(i,j\right)\in \cal C   \label{flowDemandLevel}\\
& \sum_{l \in L} W_{k}^{l-1}Y_{k}^{l} \leq &&\nonumber \\
& \qquad\sum_{\left(i,j\right) \in \cal C} \left(g_{ijkk} + \sum_{m \in K: m \neq k} \left(g_{ijkm} + g_{ijmk}\right)\right) \leq \sum_{l \in L} W_{k}^{l} Y_{k}^{l} \qquad && k \in K \label{YandgBound}\\
& \sum_{k,m\in K} T_{ijkm} x_{ijkm} + \sum_{k \in K}\sum_{l \in L} h_{k}^{l} \left(f_{ijk}^{l} + s_{ijk}^{l} + o_{ijk}^{l} \right) \leq \sum_{r\in R}H_{ij}^{r} \beta_{ij}^{r}\qquad &&  \left(i,j\right)\in {\cal C} \label{totalServiceTimeBounds}\\
& \sum_{\substack{k\in K:\\ k\ne i}}\sum_{m\in K}x_{ijkm}+\sum_{l\in L}Y_i^l \leq 1 &&  (i,j)\in\mathcal{C}\label{CA1_1}\\
& \sum_{\substack{k\in K}}\sum_{\substack{m\in K:\\m\ne j}}x_{ijkm}+\sum_{l\in L}Y_j^l \leq 1 &&  (i,j)\in\mathcal{C}\label{CA2_1}\\
&  g_{ijkm} \geq 0 \quad && \left(i,j\right) \in {\cal C}, k,m \in K \label{1domaing}\\
&  x_{ijkm} \in \left\{0,1\right\} \quad && \left(i,j\right) \in {\cal C}, k,m \in K \label{1domainX}\\
&  o_{ijk}^l, f_{ijk}^l, s_{ijk}^l \in \left\{0,1\right\} \quad && \left(i,j\right) \in {\cal C}, l \in L, k \in K \label{1domain-binary}\\
&  Y_{k}^{l} \in \left\{0,1\right\} \quad && k \in K, l \in L  \label{1domainY}\\
&  \beta_{ij}^{r} \in \left\{0,1\right\} \quad && \left(i,j\right) \in {\cal C}, r \in R.  \label{1domainBeta}
\end{align}
\end{small}

 Constraints  \eqref{1oneFlowLevel} impose that each hub can be activated in at most one service level, whereas Constraints \eqref{demandlevels} impose that each commodity can be served at most on one of its demand levels. By Constraints \eqref{demandLevelPath}, when commodity $(i,j)\in {\cal C}$ is served, then it is allocated to exactly one routing path. The role of constraints \eqref{yandxBound}-\eqref{yhatandxBound} is to identify the (\textit{first}, \textit{second} or \textit{third}) legs of the routing paths. Constraints \eqref{xAndzBound} are an adaptation of those proposed by \cite{MARIN2006274}, which relate the arcs traversed by routing paths with open hubs. Constraints \eqref{yAndYBound} play a similar role, by relating the endnodes of each of the tree legs with open hubs; note that we can now specify the service level of the involved hubs.
 Constraints  \eqref{gxbound} guarantee that the flow corresponding to each commodity circulates only through arcs of the corresponding routing paths and \eqref{flowDemandLevel} establish the specific amount of flow corresponding to each served commodity. Inequalities \eqref{YandgBound} guarantee that the flows entering the activated hubs, satisfy the capacity limits corresponding to their service levels whereas \eqref{totalServiceTimeBounds} ensure that the time limits for service routes are respected according to the  demand level of the commodities. We only impose the maximum service time constraints since it seems clear that no user would object if its service time is smaller than what is scheduled. Indeed, these constraints can be easily extended to enforce lower bounds as well.

Constraints \eqref{CA1_1}-\eqref{CA2_1} are needed to guarantee that, if the origin (destination) of a commodity $(i, j)\in\cal C$ is a hub, then the commodity is not routed through a \textit{first} (\textit{second}) hub different from $i$ ($j$). Note that this type of constraints is not needed in HLPs with no \textit{capacity-like} constraints. Indeed, if node $i$ is activated as a hub, the routing cost of $x_{ijim}$ will be smaller than that of $x_{ijkm}$ for all $k\ne i$, so no optimal solution will use the later unless there is some constraint restricting this possibility. However, in the presence of capacity constraints, the above argument no longer holds as activating $x_{ijkm}$ could help satisfying some capacity constraint involving hub $i$. This is illustrated in Example \ref{example}.

The remaining constraints determine the domain of the decision variables. Note also that, as already mentioned, the binary condition on the routing variables $x$ cannot be relaxed. The number of constraints of F1 is $2|K|+|\mathcal{C}|\left(6 + 4|K|+ |K|\times |L|+ |K|^2\right)$.

\begin{example}\label{example}
Consider the following example on the graph shown in Figure~\ref{table:example:param:costs} where the unit routing costs $c_{ij}$ are depicted next to the edges. We assume that hubs can be activated at two different service levels, i.e., $|L|=2$, and that each commodity has one single demand level, i.e. $|R|=1$. This allows us to drop all super-indices $r$ related to service levels.
The set of commodities is ${\cal{C}}=\{(1,2), (2, 4)\}$, both of them with demand 100, i.e. $w_{12}=w_{24}=100$ and unit revenue $q_{12}=q_{24}=5$.
The activation costs for hubs located at nodes 1 and 4 are $G_1^1=G_4^1=150$ and $G_1^2=G_4^2=300$, whereas the activation costs for hubs located at nodes 2 and 3 are $G_2^1=G_3^1=50$ and $G_2^2=G_3^2=150$. All hubs have capacity 100 if activated at service level 1 and capacity 200 if activated at service level 2, i.e.,  $W_k^1=100$ and $W_k^2=200$ for all $k\in \{1, 2, 3, 4\}$. All these data are summarized in Figure~\ref{table:example:param:hubs}. Let us suppose that there is no travel time limit for the served commodities and that $\alpha=0.5$.

\begin{figure}[htb]
\begin{subfigure}[b]{0.33\linewidth}
  \centering
\includegraphics[height=3cm]{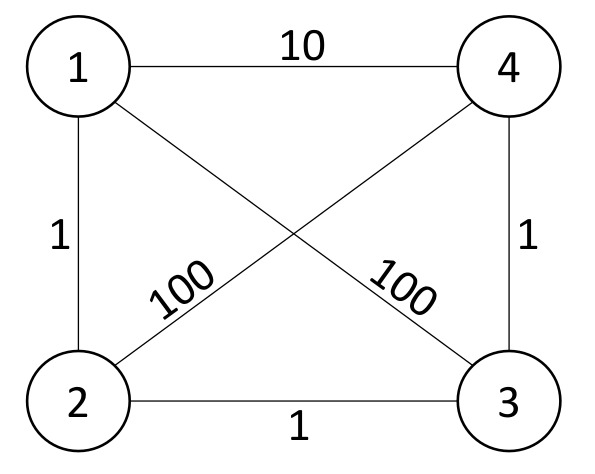}
\caption{Input graph and unit routing costs}
\label{table:example:param:costs}
\end{subfigure}
\hspace{2cm}
\begin{subfigure}[b]{0.33\linewidth}
\begin{tabular}{c|cccc}
 \multicolumn{5}{c}{}\\
\toprule
$k$ &  $W_{k}^1$ &  $W_{k}^2$ & $G_{k}^1$ & $G_{k}^2$ \\\midrule
$\{1,4\}$ & $100$ & $200$ & $150$ & $300$ \\
$\{2,3\}$ & $100$ & $200$ & $50$ & $150$\\
 \bottomrule
\end{tabular}
\caption{Data for potential hubs}
\label{table:example:param:hubs}
\end{subfigure}
\caption{Input graph and parameter values with set of commodities ${\cal{C}}=\{(1,2), (2, 4)\}$; $w_{12}=w_{24}=100$ and $q_{12}=q_{24}=5$.}
\label{table:example:param}
\end{figure}

The solution is shown in Figure~\ref{fig:example:solution}: optimal flows are depicted in Figure \ref{fig:example:sol:flows}; information on activated hubs and service levels is given in Figure \ref{fig:example:sol:hubs}; and optimal values of variables $\beta$, $x$, and $g$ and associated revenues and routing costs  in Figure \ref{fig:example:sol:com}. It is not difficult to see that any optimal solution will route both commodities so $\beta_{12}=\beta_{24}=1$.  The routing path for commodity $(1, 2)$ is just $1-2$, whereas the routing path for commodity $(2, 4)$ is $2-3-4$ (note that arc (2,4) has a very large unit routing cost), corresponding to $x_{1222}=1$ and $x_{2423}=1$, respectively. Thus, $g_{1222}=100$ and $g_{2423}=100$ so that the overall flow involved in hub 2 is $g_{1222}+g_{2423}=100+100=200$, whereas the overall flow involved in hub 3 is $g_{2423}=100$. This means that the hub located at node 2 has to be activated at level $l=2$, i.e. $Y_2^2=1$, whereas it is enough to activate at level $l=1$ the hub located at node 3, i.e. $Y_3^1=1$. The overall revenue for serving the commodities is $w_{12}q_{12}+w_{24}q_{24}=500+500=1000$; the setup costs of the activated hubs are $G_2^2+G_3^1=150+50=200$ whereas the routing cost for  $(1, 2)$ is $w_{12}\, c_{12}=100$, and the routing cost for $(2, 4)$ is $w_{24}\, \left(\alpha\, c_{23}+c_{34}\right)=100\,(0.5\cdot 1+1)=150$. In total, the objective function value of this solution is $(500+500)-200-(100+150)=1000-450=550$.

\begin{figure}[htb]
\begin{subfigure}[b]{0.33\linewidth}
\centering
\includegraphics[height=3cm]{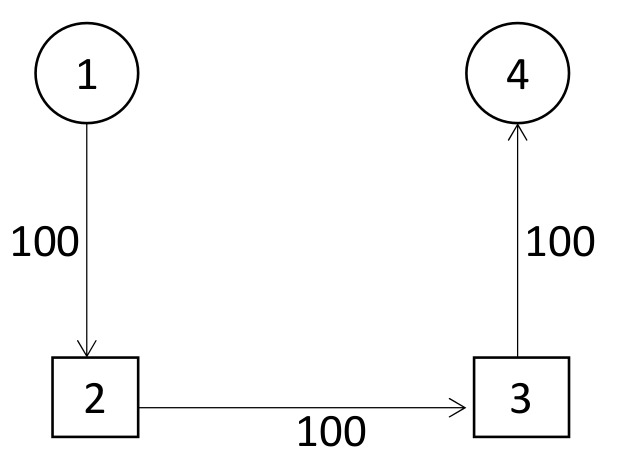}
\caption{Optimal flows}
\label{fig:example:sol:flows}
\end{subfigure}
\begin{subfigure}[b]{0.20\linewidth}
\begin{tabular}{ccc}
\toprule
$Y_k^l$  & $W_k^l$ & $G_k^l$ \\
\midrule
$Y_2^2=1$ & $200$ & $150$ \\
$Y_3^1=1$ & $100$ & $50$  \\
\bottomrule
\end{tabular}
\caption{Activated hubs}
\label{fig:example:sol:hubs}
\end{subfigure}
\hspace{0.5cm}
\begin{subfigure}[b]{0.33\linewidth}
\begin{tabular}{ccccc}
\toprule
$\beta_{ij}$ & $x_{ijkm}$  & $g_{ijkm}$ & $q_{ij}\,w_{ij}$ & $C_{ijkm}$ \\
\midrule
$\beta_{12}=1$ & $x_{1222}=1$ & $100$ & $500$ & $100$ \\
$\beta_{24}=1$ & $x_{2423}=1$ & $100$ & $500$ & $150$ \\
\bottomrule
\end{tabular}
\caption{Served commodities}
\label{fig:example:sol:com}
\end{subfigure}
\caption{Optimal solution for the example of Figure \ref{table:example:param}.}
\label{fig:example:solution}
\end{figure}

Observe that if Constraints \eqref{CA1_1}-\eqref{CA2_1} were not included, the above formulation would produce a \textit{solution} where the two commodities are still served, but the routing variables would be $x_{1222}=1$ and $x_{2433}=1$, corresponding to flow variables $g_{1222}=100$ and $g_{2433}=100$. Note that with these flow values, the overall flow involved in hub 2 would be only $g_{1222}=100$ so the hub located at node 2 could be activated at level $l=1$, i.e. $Y_1^2=1$, incurring a setup cost of $50$ (instead of the $150$ of the previous solution). No discount factor could be applied any longer to the routing cost through arc $(2, 3)$, so the routing cost of commodity $(2,4)$ would  be $100(1+1)=200$. Still, this \textit{solution} would outperform the previous one, as its  overall objective function value would be $(500+500)-(50+50)-(100+200)=1000-400=600$.

The above example illustrates the necessity of including Constraints  \eqref{CA1_1}-\eqref{CA2_1} in order to obtain a valid formulation for the HLCTDP.
\end{example}

\subsection{Formulation F2: 5-index routing variables}

In the second formulation that we develop for the HLCTDP  we use  5-index routing variables, where the new index refers to the demand level at which the concerned commodity is served.
In particular we use the following set of routing variables:

\begin{itemize}
	\item $x_{ijkm}^{r}$: binary variable, which takes the value 1, if and only if commodity $\left(i,j\right) \in \cal C$  is served at demand level $r\in R$, and it is routed first through hub $k \in K$ and then through hub $m\in K$.
\end{itemize}
There are now $\mathcal{O}\left(|\mathcal{C}|\times|K|^2\times|R|\right)$  variables $x$, but use no additional binary variables. Still, we need $\mathcal{O}\left(|\mathcal{C}|\times|K|\right)$ additional continuous variables to determine the transit times of the routing paths at the hubs that they traverse. In particular:

\begin{itemize}
	\item $\overline t_{ijk}$: continuous variable, representing the transit time of commodity $\left(i,j\right) \in \cal C$  at hub $k \in K$.
\end{itemize}

\medskip

The formulation is the following:

\begin{small}
\begin{align}
F2 \quad  \text{max } & \sum_{\left(i,j\right) \in \cal C}\sum_{r \in R}\sum_{k,m\in K} \overline C_{ijkm}^{r}x_{ijkm}^{r} - \sum_{\substack{k \in K, l \in L}} G_{k}^{l} Y_{k}^{l}&&  \label{2ofOrig}\\
\text{s.t.} &  \sum_{l\in L} Y_{k}^l \le 1 \quad && k\in K   \label{2oneFlowLevel}\\
\qquad & \sum_{r\in R} \beta_{ij}^{r} \leq 1 \quad &&  \left(i,j\right)\in {\cal C}   \label{2demandlevels1}\\
&  \beta_{ij}^{r} = \sum_{k\in K}\sum_{m\in K}x_{ijkm}^{r} \quad &&\left(i,j\right)\in {\cal C}, r\in R   \label{2beta-x}\\
& 	\sum_{r\in R}\left(x_{ijkk}^r + \sum_{\substack{m\in K:\\ m\ne k}} \left(x_{ijkm}^r + x_{ijmk}^r\right)\right) \leq
\sum_{l \in L}Y_{k}^{l}  \qquad &&  \left(i,j\right)\in {\cal C}, k \in K   \label{F3xAndzBound}\\
& \sum_{l\in L:l>1} (W_{k}^{l-1} +1) Y_{k}^{l} \leq \nonumber \\
& \qquad \sum_{\left(i,j\right) \in \cal C}\sum_{r\in R}w_{ij}^r\left(x_{ijkk}^r+\sum_{\substack{m\in K:\\ m\ne k}} (x_{ijkm}^{r}+x_{ijmk}^r)\right) \leq \sum_{l\in L} W_{k}^{l} Y_{k}^{l} \quad && k \in K \label{2YandgBound1} \\
&  \sum_{r \in R}\sum_{k,m\in K}T_{ijkm} x_{ijkm}^{r}  + \sum_{k \in K} \overline t_{ijk} \leq \sum_{r\in R}H_{ij}^{r} \beta_{ij}^{r}\quad && \left(i,j\right)\in {\cal C}   \label{2totalServiceTimeBounds1}\\
& \overline t_{ijk}\geq \sum_{l\in L}h_{k}^lY_k^l-h_k^{|L|}\left(1-\sum_{r \in R}\left(x_{ijkk}+\sum_{\substack{m\in K:\\ m\ne k}}(x_{ijkm}^r+x_{ijmk}^r)\right)\right) \quad && \left(i,j\right)\in {\cal C}, k\in K  \label{2time-ijk-LB-hub2_bis}\\
& \overline t_{ijk}\le \sum_{l \in L}h_{k}^{l} Y_{k}^{l} \quad && \left(i,j\right)\in {\cal C}, k\in K  \label{2time-ijk-UB}\\
& \sum_{\substack{k\in K:\\ k\ne i}}\sum_{m\in K}\sum_{r \in R}x_{ijkm}^{r}+\sum_{l\in L}Y_i^l \leq 1 &&  (i,j)\in\mathcal{C}\label{CA1}\\
& \sum_{k\in K}\sum_{\substack{m\in K:\\ m\ne j}}\sum_{r \in R}x_{ijkm}^{r}+\sum_{l\in L}Y_j^l \leq 1 &&  (i,j)\in\mathcal{C}\label{CA2}\\
&  \overline t_{ijk} \geq 0 \quad && \left(i,j\right) \in {\cal C}, k \in K  \label{2domainH}\\
&  x_{ijkm}^r \in \left\{0,1\right\} \quad && \left(i,j\right) \in {\cal C}, r \in R, k,m \in K \label{2domainX}\\
&  Y_{k}^{l} \in \left\{0,1\right\} \quad && k \in K, l \in L  \label{2domainY}\\
&  \beta_{ij}^{r} \in \left\{0,1\right\} \quad && \left(i,j\right) \in {\cal C}, r \in R  \label{2domainBeta}
\end{align}
\end{small}

The objective function computes the difference between the net profit associated with served commodities at their corresponding demand levels, and the setup costs of activated hubs.\\
Constraints \eqref{2oneFlowLevel} and  \eqref{2demandlevels1}  are exactly the same as Constraints \eqref{1oneFlowLevel} and \eqref{demandlevels} of F1, and guarantee that each hub can be activated in at most one service level and that each commodity can be served at most on one of its demand levels, respectively.
Constraints \eqref{2beta-x}-\eqref{F3xAndzBound} impose that a route is determined for each served commodity and that all intermediate nodes traversed by the routes are activated as hubs, respectively. Constraints \eqref{2YandgBound1} ensure that the maximum flow that circulates through each activated hub respects the  capacity associated to its service level. Constraints \eqref{2totalServiceTimeBounds1} guarantee that the routes of the served commodities respect the service time limit associated with their demand level. For this, the transit time $\overline t_{ijk}$ consumed at each traversed hub $k$, which depends on the service level at which the hub is activated, is computed in Constraints \eqref{2time-ijk-LB-hub2_bis}-\eqref{2time-ijk-UB}.
Similarly to Constraints \eqref{CA1_1} and \eqref{CA2_1}, the role of Constraints \eqref{CA1} and \eqref{CA2} is to guarantee that the commodities with origin or destination activated as a hub are properly routed. The number of constraints of F2 is $2|K|+|\mathcal{C}|\left(4 + 3 |K| + |R|\right)$.

Due to the considered objective function, the following inequalities will be satisfied by any optimal solution. Still they may reinforce the Linear Programming relaxation of $F2$ so they can be added to enhance its performance:
\begin{small}
\begin{align}
& \overline t_{ijk}\le \left(\max_{l \in L}h_{k}^{l}\right) \sum_{r\in R} \left(x_{ijkk}^{r}+\sum_{\substack{m\in K:\\ m\ne k}}x_{ijkm}^{r}+x_{ijmk}^{r}\right) \quad & \left(i,j\right)\in {\cal C}, k\in K.
\end{align}
\end{small}

\section{Optimality conditions}\label{optimalityconditions}

In this section we present some conditions that must be satisfied by optimal solutions. We use them to fix the value of some variables and exclude some non-optimal solutions.

\begin{enumerate}
\item[C1] Optimal conditions related to routing costs:
\begin{enumerate}
\item For all $i\in V$ and $(k, m)\in A$ such that $c_{i m}\le c_{i k}+\alpha\, c_{km}$, there is an optimal solution with $x_{ijkm}^r=0$ $\forall\,r\in R$. By the symmetry of the costs we also have that $x_{j i m k}^r=0$ $\forall\,r\in R$.

\item For every commodity $(i, j)\in\mathcal{C}$ and pair of hubs $(k, m)\in A$ such that $c_{ik}+c_{jm}\le c_{im}+c_{jk}$, there is an optimal solution with $x_{ijmk}^r=0$ $\forall\,r\in R$, otherwise $x_{ijkm}^r=0$ $\forall\,r\in R$.
\end{enumerate}

These conditions are adaptations to the HLCTDP of elimination criteria proposed in \citep{Boland2004}, which have been used by numerous authors.

\item[C2] Optimal conditions related to net profits: \\
If commodity $(i, j)\in\mathcal{C}$ served at demand level $r\in R$ and routed through interhub arc $(k, m)\in A$ produces no positive net profit, i.e. $q_{ij}^r-C_{ijkm}\leq 0$, then there is an optimal solution with $x_{ijkm}^r=0$.

A similar condition applies if the commodity is routed through $k$ as only hub. For every $(i, j)\in\mathcal{C}$, $k\in K$ and $r\in R$ such that $q_{ij}^r-C_{ijkk}\leq 0$, then there is an optimal solution with $x_{ijkk}^r=0$.

\item[C3] Optimal conditions related to the overall travel time constraint \eqref{2totalServiceTimeBounds1}:  \\
If commodity $(i, j)\in\mathcal{C}$ served at demand level $r\in R$ and routed through interhub arc $(k, m)\in A$ does not satisfy the time limit constraint, that is $T_{ijkm}+h_k^{min}+h_m^{min} > H_{ij}^r$, then there is an optimal solution with $x_{ijkm}^r=0$. We use the notation $h_{k}^{min} = \min_{l \in L}\left\{ h_{k}^{l}\right\}$ for each $k\in K$.

Similarly, if the commodity is routed through $k$ as only hub, when $T_{ijkk}+h_k^{min} > H_{ij}^r$ then there is an optimal solution with $x_{ijkk}^r=0$.

\end{enumerate}

Furthermore, the benchmark instances we will use in the computational experiments  (see Section \ref{instancegeneration}) satisfy the following assumptions, which would lead to further variable fixings if they were not satisfied:
\begin{enumerate}
\item[A1] $w_{ij}^r>0$,  for all $(i, j)\in\mathcal{C}$.\\
Otherwise, the commodity $(i, j)$ cannot be served at demand level $r$, thus $\beta_{ij}^r=0$ and all associated routing variables will be zero as well ($x_{ijkm}^r=0$ for all $(k, m)\in A$ and $x_{ijkk}^r=0$ for all $k\in K$).
\item[A2] $w_{ij}^r \leq W_k^{|L|}$, for all $(i, j)\in\mathcal{C}$, $r\in R$, for a given $k\in K$.\\
Otherwise, the commodity $(i,j)$ cannot be routed through hub $k$ at demand level $r$, so $x_{ijkk}^r=0$, $x_{ijkm}^r=0$ for all $(k, m)\in A$, $x_{ijmk}^r=0$ for all $(m, k)\in A$.
\item[A3] $w_{ij}^r \leq W_k^{|L|}$, for all $(i, j)\in\mathcal{C}$, $r\in R$, for all $k\in K$.\\
Otherwise, the commodity $(i,j)$ cannot be served at demand level $r$ so $\beta_{ij}^r=0$ and all associated routing variables will be zero as well.
\end{enumerate}

\section{Computational experiments}\label{experimentsresults}
In this section we describe the computational experiments we have carried out, and present and analyze the obtained numerical results. All the experiments have been run on a 2 CPUs Intel Xeon E5-2670 with 128 GB of RAM memory. All the codes are programmed in C++ using the CPLEX 22.1 library. The section is structured in several parts. First, we present the benchmark instances that we have used. Then we summarize and analyze our results. Finally, we derive some managerial insights.

\subsection{Benchmark instances} \label{instancegeneration}

Since we are not aware of any existing benchmark instances for the HLCTDP, we have obtained our instances from the well-known CAB data set\footnote{\url{https://www.researchgate.net/publication/269396247_cab100_mok}}, generating or adapting additional data when needed.
Preliminary testing has been carried out to produce data sets with optimal solutions offering a suitable balance among the different elements determining the instances.

We consider a total of $54$ instances; each instance corresponds to a specific combination of parameter values, ($\alpha$, $n$, $|L|$, $|R|$), representing discount factor, number of nodes, number of service levels, and number of demand levels, respectively.
We have considered the following parameter values: $\alpha\in\{0.2,0.5,0.8\}$; $n\in\{30,40,50\}$; $|L|\in\{1,2\}$ and $|R|\in\{1,2,3\}$. The instances generator is described in detail below.

Each original CAB instance contains data from $100$ US cities. From each original CAB instance we have obtained nine \textit{base instances}, each of them with a number of nodes $n\in \{30,40,50\}$ and $\alpha\in\{0.2,0.5,0.8\}$, selecting the $n$ nodes with largest total demand in the 100-node instance. These base instances consider one single demand level per commodity and one single service level per potential hub. In each base instance, every node is a potential hub, i.e., $|K|=n$, and the number of commodities is $|\mathcal C|=n(n-1)$. The unit routing distances of base instances, $d_{ij}$, are the ones provided by the original CAB instance, and travel costs as well as travel times are taken to have the same values as distances, i.e., $c_{ij}=t_{ij}=d_{ij}$.

To keep consistency with the notation used in the previous sections, data from base instances that must be \textit{adapted} to represent multiple demand levels or multiple service levels, are denoted using the $\,\,\widetilde{}$\quad symbol.
In particular, in a base instance the (single-level) demand and revenue of each commodity are denoted by $\widetilde w_{ij}$ and $\widetilde q_{ij}$, respectively, whereas its maximum service time is denoted by $\widetilde H_{ij}$. Similarly, the (single-level) facility setup costs of base instances are denoted by $\widetilde G_k$,  their transit times by $\widetilde h_k$, and the maximum flow through each activated hub by $\widetilde W_k$.

The specific values of the above parameters related to the commodities are the following. The demands $\widetilde w_{ij}$ are the ones provided by the original CAB sets. Since these data sets do not provide revenues, we have followed \cite{Alibeyg2016,Alibeyg2018} for generating the values $\widetilde q_{ij}$ for serving the commodities. Specifically, $\widetilde q_{ij}=\dfrac{\varphi_{ij}\,\beta_q}{|A|}\,\sum_{k,m\in K}C_{ijkm}$, where $\varphi_{ij} \sim U [0.25, 0.35]$, $\beta_q=2.5$ and $C_{ijkm}= c_{ik} + \alpha\,c_{km} + c_{mj}$. The maximum service time for each commodity is defined as $\widetilde H_{ij}=\dfrac{\beta_H}{|A|}\,\sum_{k,m\in K}T_{ijkm}$, with $\beta_H=0.4$, $T_{ijkm}=t_{ik}+\gamma\, t_{km}+t_{mj}$ and $\gamma=\alpha$.

As for the parameter values related to hubs we have used the following. The facilities setup costs $\widetilde G_k$, $k\in K$, have been obtained by multiplying those generated in \citep{DeCamargo2008} by a factor $\beta_G=3000$ \citep[see also][]{Alibeyg2016,Alibeyg2018}. The transit times at the hubs have been generated as $\widetilde h_k=\dfrac{\beta_h}{|K||\mathcal{C}|}\,\sum_{k\in K}\sum_{(i, j)\in\mathcal{C}}\left(t_{ik}+t_{kj}\right)$, with $\beta_h=0.04$, and the maximum flow through each activated hub is set as $\widetilde W_k= \beta_{W}\,\sum_{(i, j)\in\mathcal{C}}\widetilde w_{ij}$ with $\beta_{W}=0.15$.

In total we have generated nine base instances, one for each combination of the parameters $\alpha\in\{0.2,0.5,0.8\}$ and $n\in\{30,40,50\}$. Then, from each base instance we have
 obtained several \textit{full instances}, by incorporating different demand and service levels and associated data. Specifically, we have generated instances with up to three demand levels $|R|\in\{1, 2, 3\}$, and with up to two capacity levels, $|L|\in\{1, 2\}$. That is, each base instance has produced $|R|\times|L|=6$ full instances. In total, as already mentioned, $54$ instances have been generated.

The parameters related to time-sensitive demand ($w_{ij}^r$, $q_{ij}^r$, $H_{ij}^r$, for all $(i, j)\in \mathcal{C}$, $r\in R$)  are generated proportional to the values of the corresponding base instance. For instance, for every $(i, j)\in \mathcal{C}$, $w_{ij}^r = \delta^r\widetilde w_{ij}$, $\forall\,r\in R$. The factors $\delta^r$, $r\in R$, take the values shown in Table~\ref{table:parameters}.
The different demand levels are referred to as \textit{medium}, \textit{high} and \textit{low} demand levels, reflecting the type of unit revenue associated with each of them. While the single-level demand when $|R|=1$ is a \textit{medium} demand, the two demand levels when $|R|=2$  are \textit{medium} and \textit{high}, and the three demand levels when $|R|=3$  are obviously \textit{low}, \textit{medium} and \textit{high}. Recall that ``high'' revenues correspond to ``low'' maximum service times and \textit{vice-versa}.

The parameters related to hub congestion ($W_k^l$, $G_k^l$, $h_k^l$, for all $k\in K$, $l\in L$), are also generated proportional to the values of their corresponding base instances. The different service levels are referred to as \textit{medium} and \textit{high}, where ``medium'' and ``high'' now refer to the maximum flow circulating through hub $k$ (capacity of hub $k$). The single level capacity when $|L|=1$ is a \textit{medium} capacity (associated with \textit{medium} transit times), whereas the two service levels when $|L|=2$  are \textit{medium} and  \textit{high}.

\begin{table}[htb]
\centering{
\begin{tabular}{c|ccc}
\toprule
Level & $w_{ij}^r$ & $q_{ij}^r$ & $H_{ij}^r$ \\
\midrule
Low & $0.6$ & $0.8$ & $1.5$ \\
Med & $1.0$ & $1.0$ & $1.0$ \\
High & $0.4$ & $3.0$ & $0.5$ \\
 \bottomrule
\end{tabular}
\hspace{2cm}
\begin{tabular}{c|ccc}
\toprule
Level &  $W_{k}^l$ & $G_{k}^l$ & $h_{k}^l$ \\
\midrule
Med & $1.0$ & $1.0$ & $1.0$ \\
High & $2.0$ & $1.7$ & $1.25$ \\
 \bottomrule
\end{tabular}
\caption{$\delta$-factor values for demand levels (left) and service levels (right)}
\label{table:parameters}}
\end{table}

\subsection{Comparison between F1 and F2}

The first set of computational experiments focuses on the comparison of the empirical performance of the two proposed formulations F1 and F2 without any further enhancements. A time limit of 7200 s was set for each formulation and tested instance. Default parameters have been used with CPLEX. In particular, the percent tolerance for optimality is established as $0.001\%$.

The results of these experiments are summarized in Table~\ref{tabla:F1-F2}. The first four columns indicate the instance parameter values ($\alpha$, $n$, $|L|$, $|R|$).
Then, for each formulation, there is a block of six columns, which gives information on the number of variables of the formulation and its performance.
Columns \#~Bin and \#~Cont show the number of binary and continuous variables, respectively. The following columns give the percent optimality gap at termination (\%~Gap), the  percent deviation of the best solution found from the best-known solution (\%~Dev), the computing time in seconds (Time), and the number of explored nodes (Nodes). For each instance, the best-known solution is the best solution obtained over all the computational experiments (which, as we will see, corresponds to Formulation F2 enhanced with the optimality conditions). In the following, we use  ``TL'' to indicate that the solver reached time limit of 7200 s, and the symbol ``$-$'' in columns \%~Gap and \%~Dev to indicate that no feasible solution was found within the time limit. Note that in some cases the time limit was reached even if the corresponding entry in column \%~Gap shows a value 0.00. For these instances after 7200 we had $0.001\%\le \%~Gap< 0.005\%$ so the optimality of the best found solution could not be stated. To differentiate these cases from those in which the optimality of the best solution found was proven, the later are marked in bold.

Note that the number of continuous variables only depends on the number of nodes of the instances ($n$), whereas the number of binary variables also increases with the number of demand levels ($|R|$) and service levels ($|L|$), both in F1 and F2. Therefore, the difficulty for solving both formulations increases with $n$, $|R|$, and $|L|$. Even if  both formulations produced a feasible solution for all $30$ node instances, no feasible solution was obtained with the tested formulations for some $40$ and $50$ node instances.

Out of the $54$ considered instances, F1 produced a feasible solution for $35$ of them, and proved the optimality within the time limit for 8 instances. The average optimality gap at termination is $8.69\%$ with a standard deviation of $18.45\%$. The average deviation of the solutions obtained with F1 from best-known solutions is $2.71\%$ with a standard deviation of $6.41\%$.

F2 clearly outperformed F1. It produced a feasible solution for $41$ instances, and proved the optimality of the best solution found for 26 of them. The average optimality gap at termination of the $41$ instances for which a feasible solution was found is $1.53\%$, with a standard deviation of $4.42\%$. The average deviation of the solutions obtained with F2 from best-known solutions is $0.42\%$, with a standard deviation of $1.49\%$. Note that computing times are always better with F2 than with F1.

Taking into account these results for the remaining experiments we only considered formulation F2.

\begin{table}[hptb]\begin{scriptsize}\begin{center}
\caption{Comparison between F1 and F2}\label{tabla:F1-F2}
\begin{tabular}{cccc|rrrrrr|rrrrrr}
\toprule
\multicolumn{ 4}{c|}{}& \multicolumn{ 6}{c|}{F1} & \multicolumn{ 6}{c}{F2} \\ \midrule
$\alpha$ & $n$ & $|L|$ & $|R|$
& \#\text{Bin} & \#\text{Cont}& \%\text{Gap} & \%\text{Dev} & \multicolumn{1}{c}{\text{Time}} & \text{Nodes}
& \#\text{Bin} & \#\text{Cont}& \%\text{Gap} & \%\text{Dev} & \multicolumn{1}{c}{\text{Time}} & \text{Nodes} \\ \midrule
\multirow{6}{*}{0.2} 	& \multirow{6}{*}{30}   & \multirow{3}{*}{1} 	& 1 &  786510 &  26100 &           0.02 &  0.00 	&  	\text{TL}  & 136840 &   783900 &  783000 &   \textbf{0.00} & 0.00 &       6860 & 215003    \\
                   	 	&                       &                    	& 2 &  787380 &  26100 &           0.00 &  0.00 	&  	\text{TL}  &   7886 &  1567770 &  783000 &   \textbf{0.00} & 0.00 &        357 & 1848        \\
                   	 	&                       &                    	& 3 &  788250 &  26100 &          52.67 & 19.47 	&  	\text{TL}  &      0 &  2351640 &  783000 &            0.00 & 0.00 & \text{TL}  & 14410       \\ \cmidrule(lr){3- 16}
                   	 	&                       & \multirow{3}{*}{2} 	& 1 &  789150 &  26100 &           0.00 &  0.00 	&  	\text{TL}  &  23963 &   783930 &  783000 &            0.00 & 0.00 & \text{TL}  & 260395     \\
                   	 	&                       &                    	& 2 &  790020 &  26100 &          21.46 &  5.06 	&  	\text{TL}  &    758 &  1567800 &  783000 &            3.52 & 0.06 & \text{TL}  & 3239       \\
                   	 	&                       &                    	& 3 &  790890 &  26100 &          85.98 & 28.22 	&  	\text{TL}  &      0 &  2351670 &  783000 &            6.87 & 1.52 & \text{TL}  & 914         \\ \midrule
\multirow{6}{*}{0.5} 	& \multirow{6}{*}{30}   & \multirow{3}{*}{1} 	& 1 &  786510 &  26100 &  \textbf{0.00} &  0.00 	&         1011 &  59804 &   783900 &  783000 &   \textbf{0.00} & 0.00 &         27 & 2162         \\
                   	 	&                       &                    	& 2 &  787380 &  26100 &  \textbf{0.00} &  0.00 	&          343 &   1821 &  1567770 &  783000 &   \textbf{0.00} & 0.00 &         93 & 361           \\
                   	 	&                       &                    	& 3 &  788250 &  26100 &          15.72 &  2.35 	&  	\text{TL}  &      6 &  2351640 &  783000 &   \textbf{0.00} & 0.00 &       3102 & 1293         \\ \cmidrule(lr){ 3- 16}
                   	 	&                       & \multirow{3}{*}{2} 	& 1 &  789150 &  26100 &  \textbf{0.00} &  0.00 	&          998 &   9315 &   783930 &  783000 &   \textbf{0.00} & 0.00 &        166 & 1275          \\
                   	 	&                       &                    	& 2 &  790020 &  26100 &           0.66 &  0.00 	&  	\text{TL}  &   2228 &  1567800 &  783000 &   \textbf{0.00} & 0.00 &       1598 & 5740       \\
                   	 	&                       &                    	& 3 &  790890 &  26100 &          33.49 &  6.64 	&  	\text{TL}  &      0 &  2351670 &  783000 &            4.30 & 0.23 &  \text{TL} & 1874         \\ \midrule
\multirow{6}{*}{0.8} 	& \multirow{6}{*}{30}   & \multirow{3}{*}{1} 	& 1 &  786510 &  26100 &  \textbf{0.00} &  0.00 	&         1130 &  78960 &   783900 &  783000 &   \textbf{0.00} & 0.00 &         87 & 10549        \\
                   	 	&                       &                    	& 2 &  787380 &  26100 &  \textbf{0.00} &  0.00 	&          672 &  13510 &  1567770 &  783000 &   \textbf{0.00} & 0.00 &         70 & 3881         \\
                   	 	&                       &                    	& 3 &  788250 &  26100 &           1.77 &  0.34 	&  	\text{TL}  &    296 &  2351640 &  783000 &   \textbf{0.00} & 0.00 &       3061 & 48587       \\ \cmidrule(lr){ 3- 16}
                   	 	&                       & \multirow{3}{*}{2} 	& 1 &  789150 &  26100 &  \textbf{0.00} &  0.00 	&          279 &   2439 &   783930 &  783000 &   \textbf{0.00} & 0.00 &         94 & 4068           \\
                   	 	&                       &                    	& 2 &  790020 &  26100 &           0.01 &  0.00 	&  	\text{TL}  &  90631 &  1567800 &  783000 &   \textbf{0.00} & 0.00 &        871 & 142513     \\
                   	 	&                       &                    	& 3 &  790890 &  26100 &          17.24 &  8.67 	&  	\text{TL}  &     17 &  2351670 &  783000 &   \textbf{0.00} & 0.00 &       4999 & 60688       \\ \midrule\midrule
\multirow{6}{*}{0.2} 	& \multirow{6}{*}{40}   & \multirow{3}{*}{1} 	& 1 & 2502280 &  62400 &           0.01 &  0.00 	&  	\text{TL}  &  15784 &  2497600 & 2496000 &            0.01 & 0.00 &  \text{TL} & 102431  \\
                   	 	&                       &                    	& 2 & 2503840 &  62400 &          39.96 & 17.64 	&  	\text{TL}  &      0 &  4995160 & 2496000 &            0.00 & 0.00 &  \text{TL} & 17438     \\
                   	 	&                       &                    	& 3 & 2505400 &  62400 &              - &    -  	&  	\text{TL}  &      0 &  7492720 & 2496000 &               - &    - &  \text{TL} & 0             \\ \cmidrule(lr){ 3- 16}
                   	 	&                       & \multirow{3}{*}{2} 	& 1 & 2507000 &  62400 &           7.49 &  1.90 	&  	\text{TL}  &    102 &  2497640 & 2496000 &            3.16 & 0.72 &  \text{TL} & 591       \\
                   	 	&                       &                    	& 2 & 2508560 &  62400 &              - &    -  	&  	\text{TL}  &      0 &  4995200 & 2496000 &            9.83 & 2.28 &  \text{TL} & 501          \\
                   	 	&                       &                    	& 3 & 2510120 &  62400 &              - &    -  	&  	\text{TL}  &      0 &  7492760 & 2496000 &               - &    - &  \text{TL} & 0             \\ \midrule
\multirow{6}{*}{0.5} 	& \multirow{6}{*}{40}   & \multirow{3}{*}{1} 	& 1 & 2502280 &  62400 &           0.26 &  0.01 	&  	\text{TL}  &  38478 &  2497600 & 2496000 &   \textbf{0.00} & 0.00 &       4045 & 155433  \\
                   	 	&                       &                    	& 2 & 2503840 &  62400 &           4.53 &  0.34 	&  	\text{TL}  &    182 &  4995160 & 2496000 &   \textbf{0.00} & 0.00 &       1892 & 6862      \\
                   	 	&                       &                    	& 3 & 2505400 &  62400 &              - &    -  	&  	\text{TL}  &      0 &  7492720 & 2496000 &               - &    - &  \text{TL} & 0             \\ \cmidrule(lr){ 3- 16}
                   	 	&                       & \multirow{3}{*}{2} 	& 1 & 2507000 &  62400 &  \textbf{0.00} &  0.00 	&         7005 &  12187 &  2497640 & 2496000 &   \textbf{0.00} & 0.00 &       3046 & 29388   \\
                   	 	&                       &                    	& 2 & 2508560 &  62400 &          15.16 &  3.31 	&  	\text{TL}  &      0 &  4995200 & 2496000 &            0.14 & 0.00 &  \text{TL} & 2193       \\
                   	 	&                       &                    	& 3 & 2510120 &  62400 &              - &    -  	&  	\text{TL}  &      0 &  7492760 & 2496000 &               - &    - &  \text{TL} & 0             \\ \midrule
\multirow{6}{*}{0.8} 	& \multirow{6}{*}{40}   & \multirow{3}{*}{1} 	& 1 & 2502280 &  62400 &           0.00 &  0.00 	&  	\text{TL}  &  90205 &  2497600 & 2496000 &            0.00 & 0.00 &  \text{TL} & 343984  \\
                   	 	&                       &                    	& 2 & 2503840 &  62400 &           0.02 &  0.00 	&  	\text{TL}  &   1819 &  4995160 & 2496000 &   \textbf{0.00} & 0.00 &       1085 & 242      \\
                   	 	&                       &                    	& 3 & 2505400 &  62400 &              - &    -  	&  	\text{TL}  &      0 &  7492720 & 2496000 &           24.21 & 8.56 &  \text{TL} & 0           \\ \cmidrule(lr){3- 16}
                   	 	&                       & \multirow{3}{*}{2} 	& 1 & 2507000 &  62400 &  \textbf{0.00} &  0.00 	&         1467 &    200 &  2497640 & 2496000 &   \textbf{0.00} & 0.00 &        450 & 534        \\
                   	 	&                       &                    	& 2 & 2508560 &  62400 &           4.27 &  0.30 	&  	\text{TL}  &      9 &  4995200 & 2496000 &   \textbf{0.00} & 0.00 &       3004 & 297         \\
                   	 	&                       &                    	& 3 & 2510120 &  62400 &              - &    -  	&  	\text{TL}  &      0 &  7492760 & 2496000 &            0.46 & 0.02 &  \text{TL} & 722          \\ \midrule\midrule
\multirow{6}{*}{0.2} 	& \multirow{6}{*}{50}   & \multirow{3}{*}{1} 	& 1 & 6134850 & 122500 &           0.68 &  0.01 	&  	\text{TL}  &    395 &  6127500 & 6125000 &            0.00 & 0.00 &  \text{TL} & 51426    \\
                   	 	&                       &                    	& 2 & 6137300 & 122500 &              - &    -  	&  	\text{TL}  &      0 & 12254950 & 6125000 &              -  &    - &  \text{TL} & 0           \\
                   	 	&                       &                    	& 3 & 6139750 & 122500 &              - &    -  	&  	\text{TL}  &      0 & 18382400 & 6125000 &              -  &    - &  \text{TL} & 0           \\ \cmidrule(lr){3- 16}
                   	 	&                       & \multirow{3}{*}{2} 	& 1 & 6142250 & 122500 &              - &    -  	&  	\text{TL}  &      0 &  6127550 & 6125000 &           10.28 & 3.79 &  \text{TL} & 0          \\
                   	 	&                       &                    	& 2 & 6144700 & 122500 &              - &    -  	&  	\text{TL}  &      0 & 12255000 & 6125000 &              -  &    - &  \text{TL} & 0           \\
                   	 	&                       &                    	& 3 & 6147150 & 122500 &              - &    -  	&  	\text{TL}  &      0 & 18382450 & 6125000 &              -  &    - &  \text{TL} & 0           \\ \midrule
\multirow{6}{*}{0.5} 	& \multirow{6}{*}{50}   & \multirow{3}{*}{1} 	& 1 & 6134850 & 122500 &           0.00 &  0.00 	&  	\text{TL}  &  13943 &  6127500 & 6125000 &   \textbf{0.00} & 0.00 &       5973 & 102843 \\
                   	 	&                       &                    	& 2 & 6137300 & 122500 &              - &    -  	&  	\text{TL}  &      0 & 12254950 & 6125000 &   \textbf{0.00} & 0.00 &       5548 & 4183       \\
                   	 	&                       &                    	& 3 & 6139750 & 122500 &              - &    -  	&  	\text{TL}  &      0 & 18382400 & 6125000 &               - &    - & \text{TL}  & 0           \\ \cmidrule(lr){3- 16}
                   	 	&                       & \multirow{3}{*}{2} 	& 1 & 6142250 & 122500 &           1.71 &  0.56 	&  	\text{TL}  &     74 &  6127550 & 6125000 &   \textbf{0.00} & 0.00 &       5288 & 37599     \\
                   	 	&                       &                    	& 2 & 6144700 & 122500 &              - &    -  	&  	\text{TL}  &      0 & 12255000 & 6125000 &               - &    - & \text{TL}  & 0           \\
                   	 	&                       &                    	& 3 & 6147150 & 122500 &              - &    -  	&  	\text{TL}  &      0 & 18382450 & 6125000 &               - &    - & \text{TL}  & 0           \\ \midrule
\multirow{6}{*}{0.8} 	& \multirow{6}{*}{50}   & \multirow{3}{*}{1} 	& 1 & 6134850 & 122500 &           0.12 &  0.00 	&  	\text{TL}  &   1755 &  6127500 & 6125000 &   \textbf{0.00} & 0.00 &        581 & 7576     \\
                   	 	&                       &                    	& 2 & 6137300 & 122500 &           1.00 &  0.00 	&  	\text{TL}  &    113 & 12254950 & 6125000 &   \textbf{0.00} & 0.00 &       2859 & 20      \\
                   	 	&                       &                    	& 3 & 6139750 & 122500 &              - &    -  	&  	\text{TL}  &      0 & 18382400 & 6125000 &               - &    - & \text{TL}  & 0           \\ \cmidrule(lr){3- 16}
                   	 	&                       & \multirow{3}{*}{2} 	& 1 & 6142250 & 122500 &           0.00 &  0.00 	&  	\text{TL}  &    510 &  6127550 & 6125000 &   \textbf{0.00} & 0.00 &        668 & 540       \\
                   	 	&                       &                    	& 2 & 6144700 & 122500 &              - &    -  	&  	\text{TL}  &      0 & 12255000 & 6125000 &   \textbf{0.00} & 0.00 &       6350 & 204        \\
                   	 	&                       &                    	& 3 & 6147150 & 122500 &              - &    -  	&  	\text{TL}  &      0 & 18382450 & 6125000 &               - &    - & \text{TL}  & 0           \\ 
\bottomrule
\end{tabular}\end{center}\end{scriptsize}\end{table}                 

\subsection{Effect of optimality conditions on F2}

 We next analyze the effect  on formulation F2 of the optimality conditions described Section~\ref{optimalityconditions}, when applied as a preprocess step. Again a time limit of 7200 s was set and the percent tolerance for optimality was established as $0.001\%$. This information is summarized in Table~\ref{tabla:F3Preproc} where, as before,  the first four columns show the parameter values ($\alpha$, $n$, $|L|$, $|R|$).
 The following five columns summarize the effect of the elimination process: percentage of $x$-variables eliminated (\%E), contribution to the percentage of each of the three optimality conditions (\%C1, \%C2, and \%C3, respectively), and computing time in seconds (T$_\text{P}$) required by the reduction procedures. Finally, the last six columns compare some performance indicators without and with the elimination criteria. The reported values are percent optimality gaps (\%~Gap), total computing times in seconds -including preprocess time- (Time), and number of explored nodes (Nodes). We use the symbol $+$ as a superscript to denote values when applying the elimination tests.

As can be seen, the elimination tests have a remarkable effect on the number of routing variables that are eliminated, which is always beyond 95\%, entailing a  very modest computational burden, which, as expected, increases with the number of nodes and demand levels, but remains below 20 seconds in all cases. The most effective test is clearly $C1$. This is not surprising as this type of test has already shown its effectiveness in other hub location problems. The performance of $C1$ improves with the value of $\alpha$; accordingly, it is less effective for instances with $\alpha=0.2$. This can be appreciated by observing that the only three instances having optimality gaps above $0.5\%$ at termination correspond to $\alpha=0.2$.

Still, the effectiveness of $C2$ and $C3$ is not negligible as they jointly eliminate within 20-30\% of the routing variables. Moreover, they seem to somehow have a complementary effect: while they are equally effective for instances with $|R|\in\{1, 3\}$, $C3$ is clearly more effective than $C2$  for instances with $|R|=2$, for every value of $\alpha$.

The beneficial effect of these reductions can be appreciated on the optimality gaps obtained after applying the tests. Now, a feasible solution is found for all 54 instances and 34 such instances are solved to proven optimality within the time limit. The optimality gap at termination is below $3\%$ in all cases, with $51$ instances with an optimality gap below $0.5\%$. The average optimality gap is $0.10\%$ with a standard deviation of $0.42\%$. These improvements can also be appreciated in the computing times, although for some small instances with $\alpha=0.2$ the computing times are higher when applying the elimination tests than without using them.

\begin{table}[hptb]\begin{scriptsize}\begin{center}
\caption{Effect of optimality conditions on formulation F2}\label{tabla:F3Preproc}
\begin{tabular}{cccc|rcrcc|rrrrrr}
\toprule
$\alpha$ & $n$ & $|L|$ & $|R|$ & \% E & \% C1 & \% C2 & \% C3 & \text{T$_P$} &  \% \text{Gap} &\% \text{Gap}$^+$ & \text{Time} & \text{Time$^+$} & \text{Nodes} & \text{Nodes$^+$}\\ \midrule
\multirow{6}{*}{0.2} 	& \multirow{6}{*}{30} 	&\multirow{3}{*}{1} 	& 1 &  97.43 & 66 & 16$\,\,\,\,$ & 16 &  1 	 &    \textbf{0.00} &  \textbf{0.00} &       6860& 4579 & 215003 & 181555 		\\
                     	&                     	&                   	& 2 &  98.61 & 66 &  8$\,\,\,\,$ & 25 &  1 	 &    \textbf{0.00} &  \textbf{0.00} &        357&  444 &   1848 & 15539 		\\
                     	&                     	&                   	& 3 &  96.19 & 66 & 13$\,\,\,\,$ & 17 &  2 	 &             0.00 &          0.00  & \text{TL} &   TL &  14410 & 83822 		\\ \cmidrule(lr){ 3- 15}
                     	&                     	&\multirow{3}{*}{2} 	& 1 &  97.43 & 66 & 16$\,\,\,\,$ & 16 &  1 	 &             0.00 &          0.01  & \text{TL} &   TL & 260395 & 264976 		\\
                     	&                     	&                   	& 2 &  98.61 & 66 &  8$\,\,\,\,$ & 25 &  1 	 &             3.52 &  \textbf{0.00} & \text{TL} & 5940 & 3239   & 93937 		\\
                     	&                     	&                   	& 3 &  96.19 & 66 & 13$\,\,\,\,$ & 17 &  2 	 &             6.87 &          0.00  & \text{TL} &   TL & 914    & 15206 		\\ \midrule
\multirow{6}{*}{0.5} 	& \multirow{6}{*}{30} 	&\multirow{3}{*}{1} 	& 1 &  98.38 & 71 & 13$\,\,\,\,$ & 13 &  1 	 &    \textbf{0.00} &  \textbf{0.00} &         27&   87 & 2162   & 5527 			\\
                     	&                     	&                   	& 2 &  99.12 & 71 &  7$\,\,\,\,$ & 21 &  1 	 &    \textbf{0.00} &  \textbf{0.00} &         93&   13 & 361    & 324 			\\
                     	&                     	&                   	& 3 &  97.24 & 71 & 11$\,\,\,\,$ & 15 &  2 	 &    \textbf{0.00} &  \textbf{0.00} &       3102&   54 & 1293   & 370 			\\ \cmidrule(lr){ 3- 15}
                     	&                     	&\multirow{3}{*}{2} 	& 1 &  98.38 & 71 & 13$\,\,\,\,$ & 13 &  1 	 &    \textbf{0.00} &  \textbf{0.00} &        166&   84 & 1275   & 2595 			\\
                     	&                     	&                   	& 2 &  99.12 & 71 &  7$\,\,\,\,$ & 21 &  1 	 &    \textbf{0.00} &  \textbf{0.00} &       1598&   44 & 5740   & 614 			\\
                     	&                     	&                   	& 3 &  97.24 & 71 & 11$\,\,\,\,$ & 15 &  2 	 &             4.30 &          0.00  &  \text{TL}&   TL & 1874   & 143645 		\\ \midrule
\multirow{6}{*}{0.8} 	& \multirow{6}{*}{30} 	&\multirow{3}{*}{1} 	& 1 &  98.67 & 80 &  9$\,\,\,\,$ & 10 &  1 	 &    \textbf{0.00} &  \textbf{0.00} &         87&   81 & 10549  & 12114 			\\
                     	&                     	&                   	& 2 &  99.26 & 80 &  5$\,\,\,\,$ & 15 &  1 	 &    \textbf{0.00} &  \textbf{0.00} &         70&  130 & 3881   & 45182 		\\
                     	&                     	&                   	& 3 &  97.95 & 80 &  8$\,\,\,\,$ & 10 &  2 	 &    \textbf{0.00} &  \textbf{0.00} &       3061&  149 & 48587  & 17288 		\\ \cmidrule(lr){ 3- 15}
                     	&                     	&\multirow{3}{*}{2} 	& 1 &  98.67 & 80 &  9$\,\,\,\,$ & 10 &  1 	 &    \textbf{0.00} &  \textbf{0.00} &         94&  351 & 4068   & 49299 			\\
                     	&                     	&                   	& 2 &  99.26 & 80 &  5$\,\,\,\,$ & 15 &  1 	 &    \textbf{0.00} &  \textbf{0.00} &        871&  447 & 142513 & 134754 		\\
                     	&                     	&                   	& 3 &  97.95 & 80 &  8$\,\,\,\,$ & 10 &  2 	 &    \textbf{0.00} &  \textbf{0.00} &       4999& 1975 & 60688  & 84575 		\\ \midrule\midrule
\multirow{6}{*}{0.2} 	& \multirow{6}{*}{40} 	&\multirow{3}{*}{1} 	& 1 &  97.53 & 67 & 16$\,\,\,\,$ & 15 &  2 	 &             0.01 &          0.01  &  \text{TL}&   TL & 102431 & 100540 		\\
                     	&                     	&                   	& 2 &  98.69 & 67 &  8$\,\,\,\,$ & 24 &  4 	 &             0.00 &          0.01  &  \text{TL}&   TL & 17438  & 105179 		\\
                     	&                     	&                   	& 3 &  96.32 & 67 & 13$\,\,\,\,$ & 17 &  6 	 &                - &          0.00  &  \text{TL}&   TL & 0      & 7057 		\\ \cmidrule(lr){ 3- 15}
                     	&                     	&\multirow{3}{*}{2} 	& 1 &  97.53 & 67 & 16$\,\,\,\,$ & 15 &  2 	 &             3.16 &  \textbf{0.00} &  \text{TL}& 4024 & 591    & 9310 		\\
                     	&                     	&                   	& 2 &  98.69 & 67 &  8$\,\,\,\,$ & 24 &  4 	 &             9.83 &          0.76  &  \text{TL}&   TL & 501    & 1939 		\\
                     	&                     	&                   	& 3 &  96.32 & 67 & 13$\,\,\,\,$ & 17 &  7 	 &                - &          1.67  &  \text{TL}&   TL & 0      & 538 		\\ \midrule
\multirow{6}{*}{0.5} 	& \multirow{6}{*}{40} 	&\multirow{3}{*}{1} 	& 1 &  98.36 & 73 & 13$\,\,\,\,$ & 12 &  2 	 &    \textbf{0.00} &          0.00  &       4045&   TL & 155433 & 438893 		\\
                     	&                     	&                   	& 2 &  99.13 & 73 &  7$\,\,\,\,$ & 20 &  4 	 &    \textbf{0.00} &  \textbf{0.00} &       1892&  207 & 6862   & 1923 			\\
                     	&                     	&                   	& 3 &  97.29 & 73 & 11$\,\,\,\,$ & 14 &  6 	 &                - &          0.00  &  \text{TL}&   TL & 0      & 39482 		\\ \cmidrule(lr){ 3- 15}
                     	&                     	&\multirow{3}{*}{2} 	& 1 &  98.36 & 73 & 13$\,\,\,\,$ & 12 &  2 	 &    \textbf{0.00} &  \textbf{0.00} &       3046& 3235 & 29388  & 61120 		\\
                     	&                     	&                   	& 2 &  99.13 & 73 &  7$\,\,\,\,$ & 20 &  4 	 &             0.14 &  \textbf{0.00} &  \text{TL}&  897 & 2193   & 11638 		\\
                     	&                     	&                   	& 3 &  97.29 & 73 & 11$\,\,\,\,$ & 14 &  6 	 &                - &          0.04  &  \text{TL}&   TL & 0      & 3034 		\\ \midrule
\multirow{6}{*}{0.8} 	& \multirow{6}{*}{40} 	&\multirow{3}{*}{1} 	& 1 &  98.66 & 81 &  9$\,\,\,\,$ &  8 &  2 	 &             0.00 &          0.00  &  \text{TL}&   TL & 343984 & 710295 		\\
                     	&                     	&                   	& 2 &  99.28 & 81 &  4$\,\,\,\,$ & 14 &  4 	 &    \textbf{0.00} &  \textbf{0.00} &       1085&   39 & 242    & 200 			\\
                     	&                     	&                   	& 3 &  98.03 & 81 &  7$\,\,\,\,$ &  9 &  6 	 &            24.21 &  \textbf{0.00} &  \text{TL}&  350 & 0      & 1000 			\\ \cmidrule(lr){ 3- 15}
                     	&                     	&\multirow{3}{*}{2} 	& 1 &  98.66 & 81 &  9$\,\,\,\,$ &  8 &  2 	 &    \textbf{0.00} &  \textbf{0.00} &        450&  151 & 534    & 1127 			\\
                     	&                     	&                   	& 2 &  99.28 & 81 &  4$\,\,\,\,$ & 14 &  4 	 &    \textbf{0.00} &  \textbf{0.00} &       3004&   78 & 297    & 104 			\\
                     	&                     	&                   	& 3 &  98.03 & 81 &  7$\,\,\,\,$ &  9 &  6 	 &             0.46 &  \textbf{0.00} &  \text{TL}& 1194 & 722    & 1911 			\\ \midrule\midrule
\multirow{6}{*}{0.2} 	& \multirow{6}{*}{50} 	&\multirow{3}{*}{1} 	& 1 &  98.15 & 67 & 16$\,\,\,\,$ & 15 &  6 	 &             0.00 &          0.00  &  \text{TL}&   TL & 51426  & 51359 		\\
                     	&                     	&                   	& 2 &  99.03 & 67 &  8$\,\,\,\,$ & 24 & 11 	 &               -  &  \textbf{0.00} &  \text{TL}&  676 & 0      & 1269 		\\
                     	&                     	&                   	& 3 &  96.86 & 67 & 13$\,\,\,\,$ & 17 & 18 	 &               -  &          0.00  &  \text{TL}&   TL & 0      & 10773		\\ \cmidrule(lr){ 3- 15}
                     	&                     	&\multirow{3}{*}{2} 	& 1 &  98.15 & 67 & 16$\,\,\,\,$ & 15 &  6 	 &            10.28 &          0.19  &  \text{TL}&   TL & 0      & 560 		\\
                     	&                     	&                   	& 2 &  99.03 & 67 &  8$\,\,\,\,$ & 24 & 11 	 &               -  &          0.00  &  \text{TL}&   TL & 0      & 1150 		\\
                     	&                     	&                   	& 3 &  96.86 & 67 & 13$\,\,\,\,$ & 17 & 17 	 &               -  &          2.53  &  \text{TL}&   TL & 0      & 11 		\\ \midrule
\multirow{6}{*}{0.5} 	& \multirow{6}{*}{50} 	&\multirow{3}{*}{1} 	& 1 &  98.84 & 74 & 13$\,\,\,\,$ & 12 &  6 	 &    \textbf{0.00} &  \textbf{0.00} &       5973& 1339 & 102843 & 17330 		\\
                     	&                     	&                   	& 2 &  99.39 & 74 &  7$\,\,\,\,$ & 19 & 10 	 &    \textbf{0.00} &  \textbf{0.00} &       5548&  485 & 4183   & 9112 		\\
                     	&                     	&                   	& 3 &  97.77 & 74 & 11$\,\,\,\,$ & 14 & 17 	 &                - &          0.01  & \text{TL} &   TL & 0      & 9956 		\\ \cmidrule(lr){ 3- 15}
                     	&                     	&\multirow{3}{*}{2} 	& 1 &  98.84 & 74 & 13$\,\,\,\,$ & 12 &  6 	 &    \textbf{0.00} &  \textbf{0.00} &       5288& 2652 & 37599  & 32383 		\\
                     	&                     	&                   	& 2 &  99.39 & 74 &  7$\,\,\,\,$ & 19 & 10 	 &                - &  \textbf{0.00} & \text{TL} & 6341 & 0      & 38826		\\
                     	&                     	&                   	& 3 &  97.77 & 74 & 11$\,\,\,\,$ & 14 & 17 	 &                - &          0.16  & \text{TL} &   TL & 0      & 512 		\\ \midrule
\multirow{6}{*}{0.8} 	& \multirow{6}{*}{50} 	&\multirow{3}{*}{1} 	& 1 &  99.11 & 82 &  9$\,\,\,\,$ &  8 &  6 	 &    \textbf{0.00} &  \textbf{0.00} &        581& 2198 & 7576   & 6946 			\\
                     	&                     	&                   	& 2 &  99.52 & 82 &  4$\,\,\,\,$ & 13 & 10 	 &    \textbf{0.00} &  \textbf{0.00} &       2859&   75 & 20     & 344 		\\
                     	&                     	&                   	& 3 &  98.45 & 82 &  7$\,\,\,\,$ &  9 & 17 	 &                - &  \textbf{0.00} & \text{TL} &  567 & 0      & 300 			\\ \cmidrule(lr){ 3- 15}
                     	&                     	&\multirow{3}{*}{2} 	& 1 &  99.11 & 82 &  9$\,\,\,\,$ &  8 &  6 	 &    \textbf{0.00} &  \textbf{0.00} &        668&  241 & 540    & 630 			\\
                     	&                     	&                   	& 2 &  99.52 & 82 &  4$\,\,\,\,$ & 13 & 10 	 &    \textbf{0.00} &  \textbf{0.00} &       6350&  191 & 204    & 284 		\\
                     	&                     	&                   	& 3 &  98.45 & 82 &  7$\,\,\,\,$ &  9 & 17 	 &                - &  \textbf{0.00} & \text{TL} & 1292 & 0      & 404 		\\ 
\bottomrule
\end{tabular}\end{center}\end{scriptsize}\end{table}                 

\subsection{Managerial insights}

The purpose of this section is to derive some insight on the structure of the solutions of the HLCTDP, which can be useful in a decision making process. For this, we study some characteristics of the best solutions found, like the number of activated hubs, the number of commodities served, and the contribution to the objective function of its different terms: net profits and costs.

Table~\ref{tabla:SOL}, gives information on best-known solutions. Again, the first four columns give the parameter values ($\alpha$, $n$, $|L|$, $|R|$). An asterisk is used in the following column to indicate solutions whose optimality was not proven within the time limit. The information related to hubs is detailed in the four following columns: number of activated hubs (\# H), percentage of hubs installed at medium and high capacity (\% M  and \% H, respectively), and percentage of installed capacity which was actually used (\% Oc). The next four columns summarize information related to the commodities: percentage of served commodities (\% Served), and percentage served at each revenue level: high, medium,  and low  (\% H, \% M, and \% L, respectively). Finally, we give the objective function value: net profit (value in the entry  multiplied by a factor of $E+09$) and percentage of overall costs corresponding to travel and installation costs (\% T and \% I, respectively).

As can be seen, for a fixed number of nodes $n$, the number of activated hubs and the percentage of served commodities decrease as $\alpha$ increases. This indicates that the \textit{advantage} of the interhub network decreases as the discount factor  decreases, and activated hubs are no longer profitable. For some configurations, when $\alpha=0.8$, high capacity hubs are no longer effective and only medium capacity hubs are installed. However, in such cases the percentage of hub occupancy increases together with the net profit, indicating that the installed capacity is used more efficiently.

\begin{table}[hptb]\begin{scriptsize}\begin{center}
\caption{Some characteristics of best-known solutions. Asterisks indicate non-optimal solutions.}\label{tabla:SOL}
\begin{tabular}{ccccc|cccc|cccc|rcc}
\toprule
	& 	& 	& 	&  	& \multicolumn{ 4}{ c|}{Hubs (capacity) } & \multicolumn{ 4}{c|}{Commodities (revenue)} &\multicolumn{1}{c}{Profit }& \multicolumn{2}{c}{Cost} \\ \midrule
$\alpha$ & $n$ & $|L|$ & $|R|$ & &
\# H & \% M & \% H & \multicolumn{1}{c}{ \%Oc}&
\multicolumn{1}{|c}{\% Served} & \% H & \% M & \% L &
\multicolumn{1}{c}{E+09}& \% T & \% I \\ \midrule
\multirow{6}{*}{0.2} 	& \multirow{6}{*}{30}   & \multirow{3}{*}{1} 	& 1 &   & 7 & 100 &  - & 66 & 56 & - & 100 & - 	& 3.62 & 35 & 65 		\\
                     	&                       &                    	& 2 &   & 8 & 100 &  - & 49 & 66 & 24 & 76 & - 	& 7.14 & 29 & 71 		\\
                     	&                       &                    	& 3 & *	& 8 & 100 &  - & 53 & 92 & 17 & 56 & 27 	& 7.29  & 32 & 68 		\\ \cmidrule(lr){ 3- 16}
                     	&                       & \multirow{3}{*}{2} 	& 1 & *	& 8 &  75 & 25 & 58 & 69 & - & 100 & - 	& 4.35 & 32 & 68  \\
                     	&                       &                    	& 2 &   & 7 &  86 & 14 & 53 & 65 & 18 & 82 & - 	& 7.21 & 32 & 68 		\\
                     	&                       &                    	& 3 & *	& 7 &  86 & 14 & 57 & 86 & 13 & 62 & 25 	& 7.31 & 36 & 64 		\\ \midrule
\multirow{6}{*}{0.5} 	& \multirow{6}{*}{30}   & \multirow{3}{*}{1} 	& 1 &   & 4 & 100 &  - & 85 & 31 & - & 100 & - 	& 4.29 & 48 & 52 		\\
                     	&                       &                    	& 2 &   & 5 & 100 &  - & 57 & 42 & 28 & 72 & - 	& 6.62 & 38 & 62 		\\
                     	&                       &                    	& 3 &   & 5 & 100 &  - & 70 & 72 & 16 & 41 & 44 	& 7.03 & 49 & 51 		\\ \cmidrule(lr){ 3- 16}
                     	&                       & \multirow{3}{*}{2} 	& 1 &   & 5 &  80 & 20 & 68 & 36 & - & 100 & - 	& 4.54 & 41 & 59 		\\
                     	&                       &                    	& 2 &   & 5 & 100 &  0 & 57 & 42 & 28 & 72 & - 	& 6.62 & 38 & 62 		\\
                     	&                       &                    	& 3 & *	& 6 &  83 & 17 & 64 & 83 & 12 & 44 & 44 	& 7.04 & 47 & 53 		\\ \midrule
\multirow{6}{*}{0.8} 	& \multirow{6}{*}{30}   & \multirow{3}{*}{1} 	& 1 &   & 5 & 100 &  - & 82 & 41 & - & 100 & - 	& 5.67 & 47 & 53 		\\
                     	&                       &                    	& 2 &   & 5 & 100 &  - & 65 & 49 & 31 & 69 & - 	& 8.45 & 46 & 54 		\\
                     	&                       &                    	& 3 &   & 5 & 100 &  - & 79 & 77 & 20 & 44 & 36 	& 8.87 & 56 & 44 		\\ \cmidrule(lr){ 3- 16}
                     	&                       & \multirow{3}{*}{2} 	& 1 &   & 4 &  75 & 25 & 78 & 38 & - & 100 & - 	& 6.06 & 49 & 51 		\\
                     	&                       &                    	& 2 &   & 5 & 100 &  0 & 65 & 49 & 31 & 69 & - 	& 8.45 & 46 & 54 		\\
                     	&                       &                    	& 3 &   & 5 & 100 &  0 & 79 & 77 & 20 & 44 & 36 	& 8.87 & 56 & 44 		\\ \midrule\midrule 
\multirow{6}{*}{0.2} 	& \multirow{6}{*}{40}   & \multirow{3}{*}{1} 	& 1 & *	& 9 & 100 &  - & 61 & 70 & - & 100 & - 	& 5.07 & 39 & 61 		\\
                     	&                       &                    	& 2 & *	& 8 & 100 &  - & 49 & 67 & 18 & 82 & - 	& 8.75 & 34 & 66 		\\
                     	&                       &                    	& 3 & *	& 8 & 100 &  - & 53 & 85 & 14 & 64 & 22 	& 8.91 & 37 & 63 		\\ \cmidrule(lr){ 3- 16}
                     	&                       & \multirow{3}{*}{2} 	& 1 &   & 8 &  75 & 25 & 58 & 64 & - & 100 & - 	& 5.70 & 36 & 64 		\\
                     	&                       &                    	& 2 & *	& 8 &  88 & 13 & 50 & 73 & 13 & 87 & - 	& 8.78 & 36 & 64 		\\
                     	&                       &                    	& 3 & *	& 8 &  88 & 13 & 53 & 94 & 12 & 64 & 24 	& 8.99 & 40 & 60 		\\ \midrule
\multirow{6}{*}{0.5} 	& \multirow{6}{*}{40}   & \multirow{3}{*}{1} 	& 1 & *	& 6 & 100 &  - & 68 & 45 & - & 100 & - 	& 5.17 & 43 & 57 		\\
                     	&                       &                    	& 2 &   & 5 & 100 &  - & 59 & 42 & 19 & 81 & - 	& 7.73 & 45 & 55 		\\
                     	&                       &                    	& 3 & *	& 6 & 100 &  - & 69 & 85 & 10 & 48 & 42 	& 8.22 & 53 & 47 		\\ \cmidrule(lr){ 3- 16}
                     	&                       & \multirow{3}{*}{2} 	& 1 &   & 6 &  83 & 17 & 66 & 51 & - & 100 & - 	& 5.70 & 44 & 56 		\\
                     	&                       &                    	& 2 &   & 5 & 100 &  0 & 59 & 42 & 19 & 81 & - 	& 7.73 & 45 & 55 		\\
                     	&                       &                    	& 3 & *	& 6 &  83 & 17 & 65 & 83 & 9 & 48 & 42 	& 8.27 & 51 & 49 		\\ \midrule
\multirow{6}{*}{0.8} 	& \multirow{6}{*}{40}   & \multirow{3}{*}{1} 	& 1 & *	& 6 & 100 &  - & 73 & 49 & - & 100 & - 	& 6.92 & 49 & 51 		\\
                     	&                       &                    	& 2 &   & 5 & 100 &  - & 64 & 49 & 27 & 73 & - 	& 9.73 & 48 & 52 		\\
                     	&                       &                    	& 3 &   & 5 & 100 &  - & 78 & 77 & 17 & 46 & 37 	& 10.25 & 58 & 42 		\\ \cmidrule(lr){ 3- 16}
                     	&                       & \multirow{3}{*}{2} 	& 1 &   & 5 &  80 & 20 & 73 & 50 & - & 100 & - 	& 7.39 & 50 & 50 		\\
                     	&                       &                    	& 2 &   & 5 & 100 &  0 & 64 & 49 & 27 & 73 & - 	& 9.73 & 48 & 52 		\\
                     	&                       &                    	& 3 &   & 5 & 100 &  0 & 78 & 77 & 17 & 46 & 37 	& 10.25 & 58 & 42 		\\ \midrule\midrule
\multirow{6}{*}{0.2} 	& \multirow{6}{*}{50}   & \multirow{3}{*}{1} 	& 1 & *	& 9 & 100 &  - & 61 & 66 & - & 100 & - 	& 6.86 & 44 & 56 		\\
                     	&                       &                    	& 2 &   & 9 & 100 &  - & 47 & 70 & 14 & 86 & - 	& 10.66 & 36 & 64 		\\
                     	&                       &                    	& 3 & *	& 9 & 100 &  - & 50 & 91 & 11 & 65 & 24 	& 10.79 & 39 & 61 	\\ \cmidrule(lr){ 3- 16}
                     	&                       & \multirow{3}{*}{2} 	& 1 & *	& 8 &  75 & 25 & 57 & 60 & - & 100 & - 	& 7.45 & 39 & 61 		\\
                     	&                       &                    	& 2 & *	& 9 & 100 &  0 & 47 & 70 & 14 & 86 & - 	& 10.63 & 36 & 64 		\\
                     	&                       &                    	& 3 & *	& 9 & 100 &  0 & 50 & 91 & 11 & 65 & 24 	& 10.79 & 39 & 61 	\\ \midrule
\multirow{6}{*}{0.5} 	& \multirow{6}{*}{50}   & \multirow{3}{*}{1} 	& 1 &   & 7 & 100 &  - & 67 & 47 & - & 100 & - 	& 6.99 & 47 & 53 		\\
                     	&                       &                    	& 2 &   & 6 & 100 &  - & 62 & 48 & 17 & 83 & - 	& 9.61 & 48 & 52 		\\
                     	&                       &                    	& 3 & *	& 7 & 100 &  - & 63 & 86 & 10 & 45 & 45 	& 10.09 & 51 & 49 	\\ \cmidrule(lr){ 3- 16}
                     	&                       & \multirow{3}{*}{2} 	& 1 &   & 6 &  83 & 17 & 67 & 48 & - & 100 & - 	& 7.48 & 48 & 52 		\\
                     	&                       &                    	& 2 &   & 6 & 100 &  0 & 62 & 48 & 17 & 83 & - 	& 9.61 & 48 & 52 		\\
                     	&                       &                    	& 3 & *	& 6 &  83 & 17 & 66 & 84 & 8 & 46 & 46 	& 10.22 & 55 & 45 		\\ \midrule
\multirow{6}{*}{0.8} 	& \multirow{6}{*}{50}   & \multirow{3}{*}{1} 	& 1 &   & 6 & 100 &  - & 73 & 52 & - & 100 & - 	& 8.56 & 52 & 48 		\\
                     	&                       &                    	& 2 &   & 6 & 100 &  - & 59 & 51 & 22 & 78 & - 	& 11.50 & 48 & 52 		\\
                     	&                       &                    	& 3 &   & 6 & 100 &  - & 72 & 84 & 13 & 50 & 37 	& 12.15 & 59 & 41 	\\ \cmidrule(lr){ 3- 16}
                     	&                       & \multirow{3}{*}{2} 	& 1 &   & 6 &  83 & 17 & 67 & 52 & - & 100 & - 	& 8.93 & 50 & 50 		\\
                     	&                       &                    	& 2 &   & 6 & 100 &  0 & 59 & 51 & 22 & 78 & - 	& 11.50 & 48 & 52 		\\
                     	&                       &                    	& 3 &   & 6 & 100 &  0 & 72 & 84 & 13 & 50 & 37 	& 12.15 & 59 & 41 	\\ 
\bottomrule
\end{tabular}\end{center}\end{scriptsize}\end{table}                 

Regarding the objective function value, in general, for the same configuration of $\alpha$, $|L|$ and $|R|$, as the number of nodes in the network increases (the same network is expanded with more commodities) the net profit increases as well. Roughly speaking this indicates that, as the network size increases, overall installation plus service costs increase at a lower rate than revenues.
For most instances, hub activation costs are higher than travel costs. Still, the difference between the two types of costs reduces as the discount factor decreases (and the number of activated hubs decreases as well), and in some cases the trade-off is reversed.

The structure of the obtained solutions seems to be quite \textit{stable}, when $n$ and $\alpha$ are fixed, as the number of service levels increases, even if some changes can be detected. For instance, we can observe a slight increase of the net profit. Moreover, the number of activated hubs tends to decrease and more diversity appears in terms of the installed capacity, with some high-capacity hubs replacing some medium-capacity hubs. On the contrary, no clear difference can be observed on the percentage of served demand or on the occupancy the activated capacity.

Increasing the number of demand levels has a direct effect on the potential flexibility of the decision making process, which, in its turn, affects strongly the obtained profit, since this allows solutions where the capacity installed at the activated hubs can be used more efficiently, so their occupancy increases. On the one hand, the possibility of serving commodities with higher revenues, enables to capture more profit, even if this further constrains the overall service time limit. On the other hand, serving demand with low revenue is easier since service-time constraints become less restrictive, allowing a better use of hubs' capacities.

\begin{figure}[htbp]
\begin{center}
\includegraphics[scale=0.65]{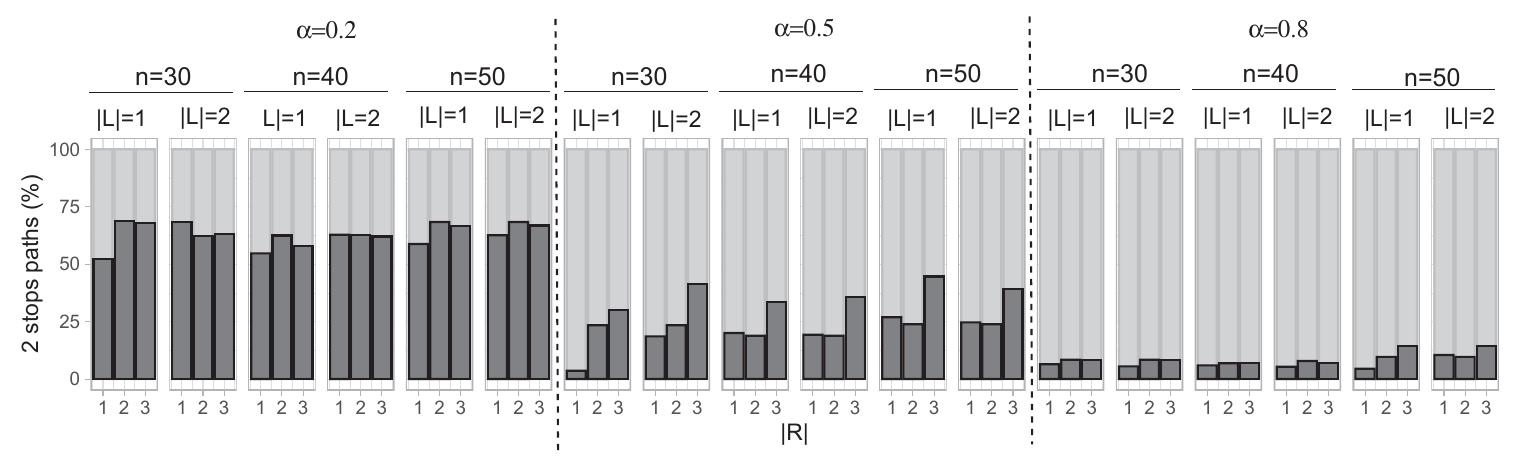}
\end{center}
\caption{Percentage of served commodities routed using two intermediate hubs.}
\label{Fig:2-stops-path}
\end{figure}

Figure~\ref{Fig:2-stops-path} shows, for each considered instance, the percentage of served commodities whose routing path use two intermediate hubs. The figure clearly shows that this percentage decreases with as the value of $\alpha$ increases. This behavior could be expected since as the discount factor is smaller (larger values of $\alpha$) interhub links become less attractive. When the value of $\alpha$ is fixed, we can observe that instances with one demand level ($|R|=1$) tend to have routing paths with fewer intermediate hubs than their counterparts with $|R|>1$. The biggest differences between the cases with $|R|=1$ and $|R|=3$ can be appreciated for the intermediate values of $\alpha=0.5$. When both $\alpha$ and $|R|$ are fixed, there seems to be a slight increase on the number of paths with two intermediate hubs for instances with $|L|=2$.

\begin{center}
\begin{figure}[htbp]
\begin{subfigure}[b]{0.33\linewidth}
  \centering
  \efbox{\includegraphics[width=0.9\linewidth]{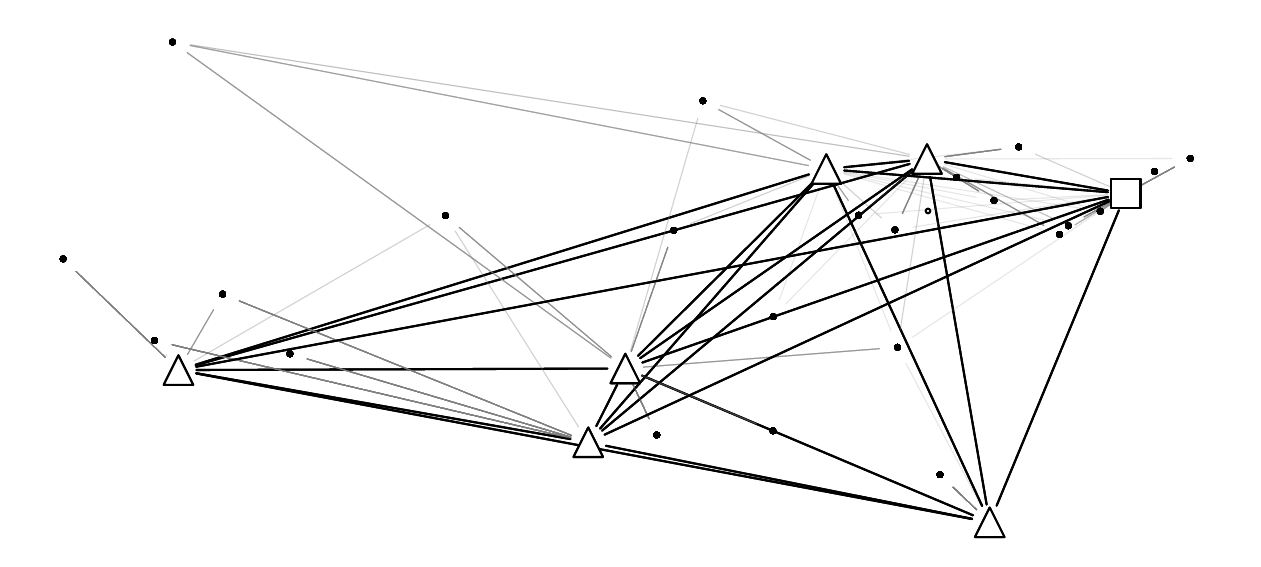} }
  \caption{$\alpha=0.20$: 7 hubs (6M + 1H)}
\end{subfigure}
\begin{subfigure}[b]{0.33\linewidth}
  \centering
  \efbox{\includegraphics[width=0.9\linewidth]{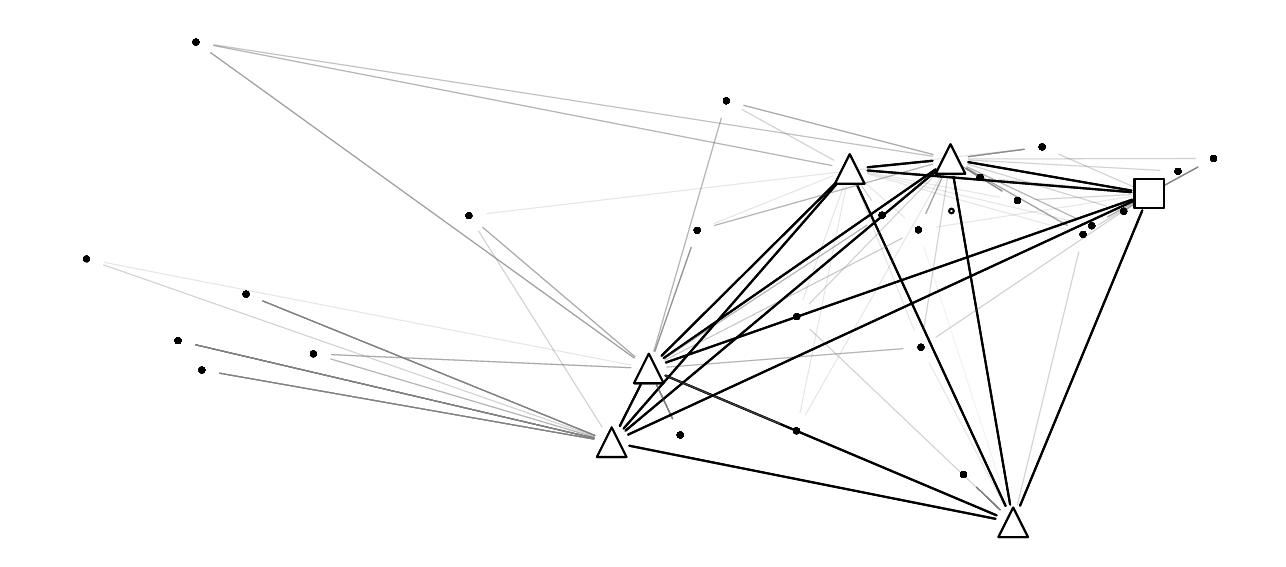}  }
  \caption{$\alpha=0.50$:  6 hubs (5M + 1H)}
\end{subfigure}
\begin{subfigure}[b]{0.33\linewidth}
  \centering
  \efbox{\includegraphics[width=0.9\linewidth]{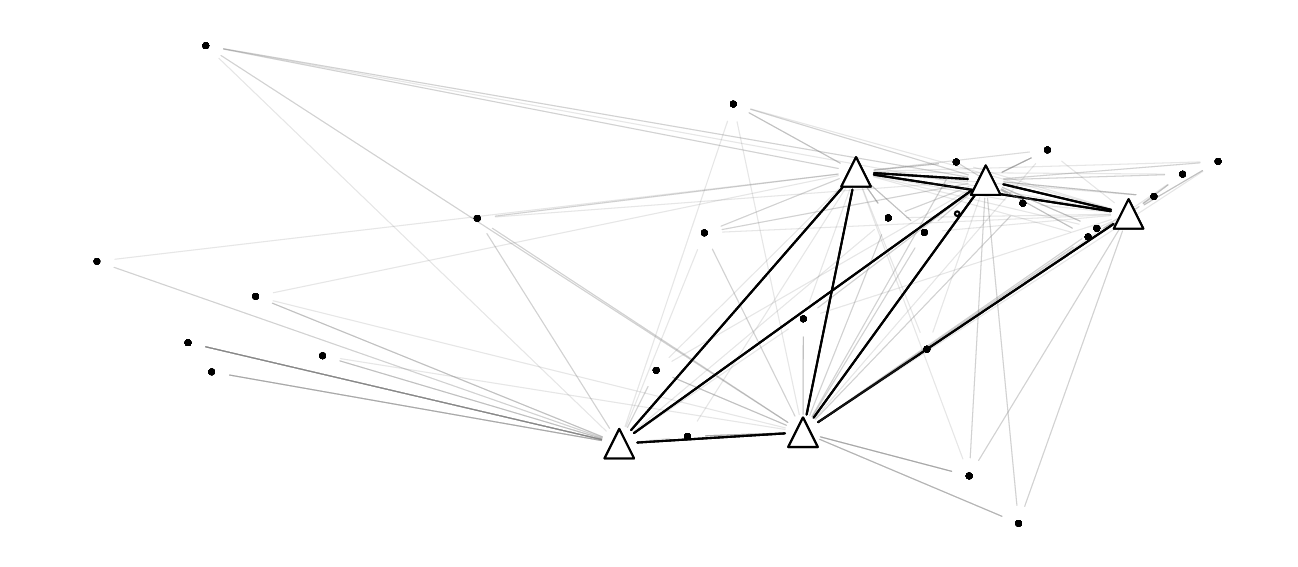} }
  \caption{$\alpha=0.80$: 5 hubs (5M + 0H)}
\end{subfigure}
\caption{Optimal solutions for $n=30$, $|L|=2$ and $|R|=3$.}
\label{Fig:Mapas}
\end{figure}
\end{center}

We finally illustrate some of the solutions characteristics by visualizing the obtained solutions for the instances with $n=30$, $|R|=3$, and $|L|=2$, for the different values of $\alpha$. These three instances were optimally solved within the maximum time limit, so in all three cases we refer to optimal solutions. Figure~\ref{Fig:Mapas} presents the optimal network for each value of $\alpha$. Activated hubs are marked with white triangles (for medium capacity) and white squares (for high capacity). Non-hub nodes are marked with small circles. Black links represent inter-hub connections whereas grey links represent  connections used as access/distribution flows. The intensity of the grey increases with the amount of flow circulating through the link. As there are no link activation costs, the backbone network of the depicted solutions is complete.

In all three cases there is only one node from which no commodity is served. It is placed in the north-east area close to the Detroit hub, which appears in all three solutions and can be identified as the medium-capacity hub which is the second left-most in the northern part of the USA.

As mentioned, the number of activated hubs decreases as $\alpha$ increases (i.e. the discount decreases). In particular, seven hubs are activated for $\alpha=0.2$ (six at medium capacity and one at high capacity), whereas the number of activated hub reduces to six for $\alpha=0.5$ (five at medium capacity and one at high capacity) and to five for $\alpha=0.8$ when all hubs are activated at medium capacity. Note that the locations for the hubs are not always the same, although four hub locations remain activated in all three solutions.

\begin{figure}[htbp]
\centering
\efbox{\includegraphics[width=0.75\linewidth]{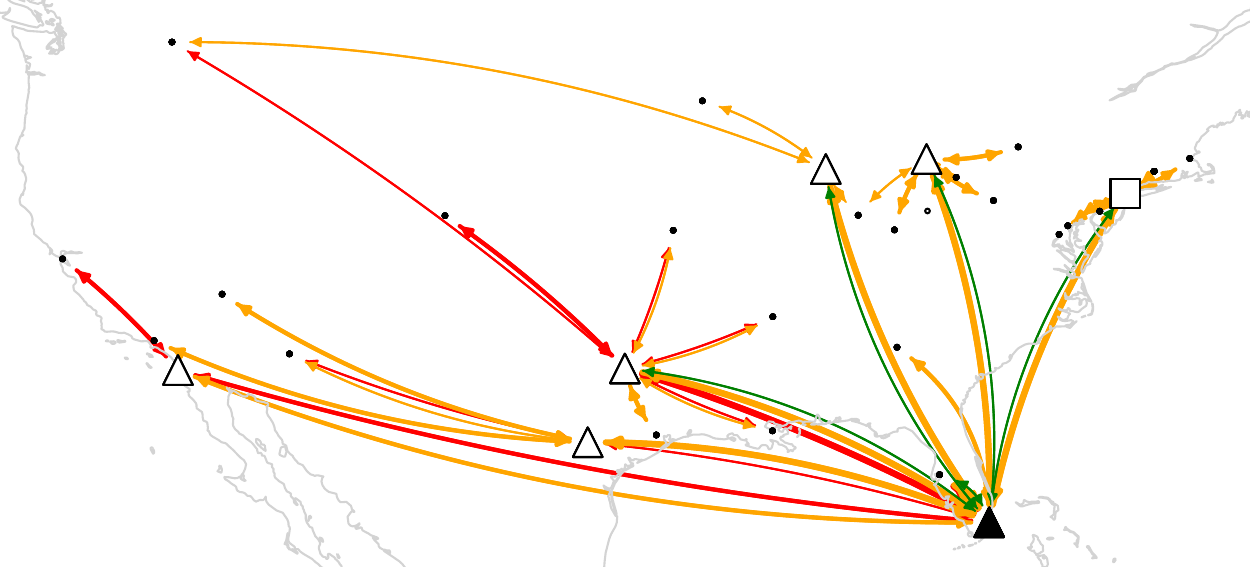} }
\caption{Flow routed through a fixed hub marked in black. (Instance with $\alpha=0.2$, $n=30$, $|L|=2$, $|R|=3$).}
\label{Fig:Mapa-Zoom}
\end{figure}

Figure~\ref{Fig:Mapa-Zoom} illustrates how flows are routed. For this we have chosen the instance with $\alpha=0.2$ and chosen a specific hub, the one located at Miami, which is marked in black, and observed the flows circulating through that specific hub. Colors indicate different demand levels: High$\,\leftrightarrow\,$green, Medium$\,\leftrightarrow\,$orange and Low$\,\leftrightarrow\,$red. The same two endnodes may be connected with links of different colors, meaning that different commodities served at different demand levels are routed using the same link. The thickness of a line increases as the flow routed through the line increases.
Green lines indicate commodities routed through two intermediate hubs at a high revenue level, so their maximum service time is \textit{low}. In this example, such connections appear only when both the origin and destination are hubs. Most of the two-hub routes corresponding to the other revenue levels (medium and low) involve at least one non-hub node.

\section{Conclusions}\label{conclusion}

Hub congestion has a very negative effect on the users of hub-and-spoke systems, since transit times at hubs directly influence overall service times from origins to destinations, which, in its turn, affect service demand. In this work we have studied how hub congestion affects time-sensitive demand for a hub location model that we call Hub Location with Congestion and Time-sensitive Demand Problem (HLCTDP).

In the HLCTDP the efficiency of a hub-and-spoke system is measured in terms of the maximum net profit it may produce. Hubs can be activated at several service levels, where each level is characterized by a maximum capacity (expressed as the amount of flow that may circulate through the hub) together with a unit transit time at the hub. In addition, we consider time-sensitive demand, by associating each commodity with several potential demand levels, each of them entailing a maximum total service time and corresponding revenues. To the best of our knowledge this is the first work where hub congestion and time-sensitive demand are jointly considered.

We have developed two alternative formulations for the HLCTDP, and several elimination tests that allow reducing the number of decision variables substantially. Both formulations use path variables to identify the routing paths of served commodities, and require a set of constraints, new in the hub location literature, that relate path variables with the hub location variables.
The new set of constraints is necessary to guarantee the consistency of the obtained solutions as, otherwise, the \textit{meaning} of the routing decision variables may be distorted due to the \textit{capacity-type} constraints derived from hub service levels and served demand levels.

Through a set of computational experiments we have analyzed the efficiency of each of the formulations as well as the effectiveness of the elimination tests. Furthermore, the results of the computational experiments have allowed us to derive interesting managerial insight on the structure of the obtained solutions and on how the different parameters affect them.

In our opinion, the HLCTDP opens several avenues for future research both from the modeling and the solution point of view. A number of extensions or variations of the HLCTDP deserve attention from the modeling perspective. Some examples are $(i)$ to consider more sophisticated costs by including link activation setup costs, or by using demand-level dependent routing costs, or $(ii)$ to consider the possibility of serving at several demand levels simultaneously. The design of efficient solution methods is a promising area of research, since off-the-shelf solvers have proven to be effective for instances with up to 40-50 nodes although, due to the large number of binary decision variables, they will no longer be able to deal with larger instances. Hence, alternative solution solution methods are needed.

\section*{Acknowledgments}

The research of Carmen-Ana Dom\'inguez and Elena Fernández has been partially funded by the Spanish Ministry of Science, Innovation and Universities  [Grant MTM-PID2019-105824GB-I00] and the Spanish Location Network [Grant MTM-RED2018-102363-T]. The research of Armin L\"{u}er-Villagra has been partially funded by the Chilean FONDECYT [Grant 1200706]. This support is gratefully acknowledged.

\bibliographystyle{elsarticle-harv}
\bibliography{References20230529-no-doi}

\end{document}